\documentclass[a4paper,12pt,twoside]{article}
\usepackage[francais]{babel}
\usepackage{amssymb, amsmath, eucal}
\usepackage[all]{xy}
\usepackage{mathrsfs}
\begin{document}
\textwidth15.5cm
\textheight22.5cm
\voffset=-13mm
\newtheorem{The}{Theorem}[section]
\newtheorem{Lem}[The]{Lemma}
\newtheorem{Prop}[The]{Proposition}
\newtheorem{Cor}[The]{Corollary}
\newtheorem{Rem}[The]{Remark}
\newtheorem{Obs}[The]{Observation}
\newtheorem{Titre}[The]{\!\!\!\! }
\newtheorem{Conj}[The]{Conjecture}
\newtheorem{Question}[The]{Question}
\newtheorem{Prob}[The]{Problem}
\newtheorem{Def}[The]{Definition}
\newtheorem{Not}[The]{Notation}
\newtheorem{Claim}[The]{Claim}
\newtheorem{Ex}[The]{Example}
\newtheorem{Com}[The]{Comment}
\newtheorem{Concl}[The]{Conclusion}
\newcommand{\C}{\mathbb{C}}
\newcommand{\R}{\mathbb{R}}
\newcommand{\N}{\mathbb{N}}
\newcommand{\Z}{\mathbb{Z}}
\newcommand{\Q}{\mathbb{Q}}

\begin{center}

{\Large\bf Limits of Projective and $\partial\bar\partial$-Manifolds under Holomorphic Deformations}

\end{center}

\begin{center}

 {\large Dan Popovici}

\end{center}

\vspace{1ex}

\noindent {\small {\bf Abstract.} We prove that if in a (smooth) holomorphic family of compact complex manifolds all the fibres, except one, are projective, then the remaining (limit) fibre must be Moishezon. In an earlier work, we proved this result under the extra assumption that the limit fibre carries a strongly Gauduchon metric. In the present paper, we remove the extra assumption by proving that if all the fibres, except one, are $\partial\bar\partial$-manifolds, then the limit fibre carries a strongly Gauduchon metric. The $\partial\bar\partial$-assumption on the generic fibre is much weaker than the projective, K\"ahler and even {\it class} ${\cal C}$ assumptions, but it implies the Hodge decomposition and symmetry, while being called the 'validity of the $\partial\bar\partial$-lemma' by many authors. Our method consists in starting off with an arbitrary smooth family $(\gamma_t)_{t\in\Delta}$ of Gauduchon metrics on the fibres $(X_t)_{t\in\Delta}$ and in correcting $\gamma_0$ in a finite number of steps to a strongly Gauduchon metric by repeated uses of the $\partial\bar\partial$-assumption on the generic fibre and of estimates of minimal $L^2$-norm solutions for $\partial$-, $\bar\partial$- and $d$-equations.}

\vspace{3ex}

\section{Introduction}

 Let $(X_t)_{t\in\Delta}$ be a holomorphic family of compact complex manifolds. This means (cf. [Kod86]) that we are given a proper holomorphic submersion $\pi: {\cal X} \rightarrow \Delta$ between complex manifolds. Thus, the fibres $X_t:=\pi^{-1}(t)$ are (smooth) compact complex manifolds of equal dimensions $n$ that vary holomorphically with $t$. Without loss of generality for our purposes in this work, we may assume that $\Delta$ is an open disc containing the origin in $\C$.
 
 The main result of this paper is the following long-conjectured statement.

\begin{The}\label{The:projdef} Let $\pi: {\cal X} \longrightarrow \Delta$ be a holomorphic family of compact complex manifolds such that the fibre $X_t:=\pi^{-1}(t)$ is projective for every $t\in\Delta^{\star}\!\!:=\Delta\setminus\{0\}$. Then $X_0:=\pi^{-1}(0)$ is Moishezon.

\end{The}

As usual, by a compact complex manifold $X$ being Moishezon we mean that there exists a projective manifold $\tilde{X}$ and a proper holomorphic bimeromorphic map (i.e. a holomorphic modification) $\mu:\tilde{X}\rightarrow X$. We know from an example of Hironaka [Hir62] that the limit fibre $X_0$ need not be projective or even K\"ahler, so Theorem \ref{The:projdef} is optimal. 

 In [Pop13], we introduced the notion of {\it strongly Gauduchon (sG) metric}: any $C^{\infty}$ positive definite $(1,\,1)$-form (i.e. any Hermitian metric) $\omega>0$ on a given compact complex $n$-dimensional manifold $X$ that satisfies the condition

\begin{equation}\partial\omega^{n-1} \hspace{1ex} \mbox{is} \hspace{1ex} \bar\partial-\mbox{exact}.\end{equation}

\noindent Unlike Gauduchon metrics (i.e. those $\omega>0$ such that $\partial\omega^{n-1}$ is $\bar\partial$-closed) which exist in the conformal class of any Hermitian metric on any compact complex manifold (cf. [Gau77]), strongly Gauduchon metrics need not exist. Compact complex manifolds that do not carry strongly Gauduchon metrics include the Calabi-Eckmann manifolds [CE53], the Tsuji manifolds [Tsu84] (see also [Pop14] for the non-sG discussion) and certain nilmanifolds and solvmanifolds (cf. e.g. [COUV16]). Compact complex manifolds that carry strongly Gauduchon metrics were termed {\it strongly Gauduchon manifolds} in [Pop13, Definition 4.1]. 

 It is under the extra (relatively weak) strongly Gauduchon assumption on the limit fibre $X_0$ that Theorem \ref{The:projdef} was proved in [Pop13].

 \begin{The}(Theorem 1.4 in [Pop13])\label{The:projdef_sG-assumption} Let $\pi: {\cal X} \longrightarrow \Delta$ be a holomorphic family of compact complex manifolds such that the fibre $X_t:=\pi^{-1}(t)$ is projective for every $t\in\Delta^{\star}$. Suppose that $X_0$ is a strongly Gauduchon manifold. Then $X_0$ is Moishezon.

\end{The}

 In the present work, we remove the strongly Gauduchon assumption on $X_0$ even in a far more general situation than the one of Theorems \ref{The:projdef} and \ref{The:projdef_sG-assumption}. The fibres $X_t$ corresponding to $t\in\Delta^\star$ will only be assumed to be {\it $\partial\bar\partial$-manifolds} in the following sense.

\begin{Def}\label{Def:dd-bar_def} A compact complex manifold $X$ is said to be a {\bf $\partial\bar\partial$-manifold} if for every bidegree $(p,\,q)$ and every $C^{\infty}$ $(p,\,q)$-form $u$ on $X$ such that $du=0$, the following exactness conditions are equivalent:

\begin{equation}\label{eqn:dd-bar_def}u\in\mbox{Im}\,d \iff u\in\mbox{Im}\,\partial \iff u\in\mbox{Im}\,\bar\partial \iff u\in\mbox{Im}\,\partial\bar\partial.\end{equation}

\end{Def}

Many authors use the term {\it cohomologically K\"ahler} for compact complex manifolds $X$ on which the implication $u\in\mbox{Im}\,d \Longrightarrow u\in\mbox{Im}\,\partial\bar\partial$ holds for all, not necessarily pure-type, forms $u$. This latter notion was first considered by Deligne, Griffiths, Morgan and Sullivan in [DGMS75] and can be seen to be equivalent to our $\partial\bar\partial$-manifold notion of Definition \ref{Def:dd-bar_def}. The classical {\it $\partial\bar\partial$-Lemma} asserts that every compact K\"ahler manifold is a $\partial\bar\partial$-manifold. The idea of studying $\partial\bar\partial$-manifolds independently of the K\"ahler realm probably goes back to [DGMS75] where any contraction of a $\partial\bar\partial$-manifold was shown to be again a $\partial\bar\partial$-manifold. In particular, {\it class} ${\cal C}$ manifolds (i.e. those manifolds that are bimeromorphically equivalent to compact K\"ahler manifolds) are $\partial\bar\partial$-manifolds. However, there exist many $\partial\bar\partial$-manifolds that are not of {\it class} ${\cal C}$ (see e.g. [Pop14, Observation 4.10], [FOU15, Theorem 5.2], [Fri17]).

On any $\partial\bar\partial$-manifold $X$, there are {\it canonical} (i.e. depending only on the complex structure of $X$) Hodge decomposition and Hodge symmetry isomorphisms in every (bi)-degree. So these manifolds behave cohomologically like compact K\"ahler ones. However, for any compact complex manifold $X$ with $\mbox{dim}_\C X\geq 3$, the following implications are strict:

\begin{eqnarray}\nonumber X \hspace{1ex} \mbox{is projective}  & \implies & X \hspace{1ex} \mbox{is K\"ahler} \hspace{1ex} \implies X \hspace{1ex} \mbox{is of class} \hspace{1ex} {\cal C} \implies \\
\nonumber &    &  X \hspace{1ex} \mbox{is a} \hspace{1ex} \partial\bar\partial-\mbox{manifold} \hspace{1ex} \implies \hspace{1ex} E_1(X) = E_\infty (X),\end{eqnarray}

\noindent where the expression $E_1(X) = E_\infty (X)$ means that the Fr\"olicher spectral sequence of $X$ degenerates at $E_1$. Moreover, every $\partial\bar\partial$-manifold is strongly Gauduchon ([Pop13]), but there are many strongly Gauduchon manifolds that are not $\partial\bar\partial$-manifolds. See e.g. [Pop14] for a rundown on these and other related facts. 

 The precise form of the result that removes the strongly Gauduchon assumption on $X_0$ from Theorem \ref{The:projdef_sG-assumption} is the following

\begin{The}\label{The:limitsG} Let $\pi : {\cal X}\rightarrow \Delta$ be a holomorphic family of compact complex manifolds such that the fibre $X_t:=\pi^{-1}(t)$ is a {\bf $\partial\bar\partial$-manifold} for every $t\in\Delta^{\star}:=\Delta\setminus\{0\}$. Then $X_0$ is a {\bf strongly Gauduchon manifold}.

\end{The}

 It is clear that the combined Theorems \ref{The:projdef_sG-assumption} and \ref{The:limitsG} prove Theorem \ref{The:projdef}. The object of this paper is thus reduced to proving Theorem \ref{The:limitsG} whose hypothesis will be supposed to hold throughout the paper. To this end, we will show that any family $(\gamma_t)_{t\in\Delta}$ of Gauduchon metrics varying in a $C^{\infty}$ way with $t$ can be modified to a family $(\rho_t)_{t\in\Delta}$ of {\it strongly Gauduchon} metrics varying in a $C^{\infty}$ way with $t$.

\vspace{2ex}

 We believe that Theorem \ref{The:limitsG} has an interest of its own besides the key role it plays in the proof of Theorem \ref{The:projdef}. It reveals a fundamental novel deformation property of $\partial\bar\partial$-manifolds and is optimal since

\vspace{1ex}

-the $\partial\bar\partial$ assumption on the fibres $X_t$ with $t\in\Delta^\star$ cannot be weakened to the strongly Gauduchon assumption on these fibres (cf. [COUV16, Theorem 5.9]), so the strongly Gauduchon property of compact complex manifolds is not closed under deformations of the complex structure. (It is, however, open -- see [Pop14].);

\vspace{1ex}

-the strongly Gauduchon conclusion on the limit fibre $X_0$ cannot be strengthened to the $\partial\bar\partial$ conclusion, so the $\partial\bar\partial$ property of compact complex manifolds is not closed under deformations (cf. [AK13] or [FOU15]).

\vspace{1ex}

 Moreover, thanks to the $\partial\bar\partial$ property being satisfied by all compact K\"ahler and all {\it class} ${\cal C}$ manifolds, Theorem \ref{The:limitsG} is expected to play a key role in future attacks on conjectural transcendental analogues of Theorem \ref{The:projdef} predicting that the projectivity assumption on the generic fibre can be weakened to the K\"ahler or the {\it class} ${\cal C}$ assumption to imply the weaker conclusion that the limit fibre $X_0$ is a {\it class} ${\cal C}$ manifold. Theorem \ref{The:limitsG} has already provided one of the key arguments to the author's proof of the deformation closedness of the Moishezon property of compact complex manifolds in [Pop10].

\section{Notation and preliminary remarks}\label{section:notation_preliminaries}

Since the base $\Delta$ of the family $(X_t)_{t\in\Delta}$ is contractible, the classical result of Ehresmann [Ehr47] ensures that all the fibres $X_t:=\pi^{-1}(t)$ are $C^{\infty}$-diffeomorphic to a fixed compact differentiable manifold $X$. So, the information contained in the family $(X_t)_{t\in\Delta}$ is equivalent to that contained in a fixed $C^\infty$ manifold $X$ equipped with a family of complex structures $(J_t)_{t\in\Delta}$ varying holomorphically with $t$. In particular, the differentiable structure $d$ of $X_t$ is independent of $t$ and for every $k$, the De Rham cohomology groups $H^k(X_t, \, \C)$ of the fibres can be identified with a fixed $\C$-vector space $H^k(X, \, \C)$. 

 However, the operators $\partial_t$ and $\bar\partial_t$, as well as the Dolbeault cohomology groups $H^{p, \, q}(X_t, \, \C)$ of the fibres $X_t$ depend on $t\in\Delta$. So the splitting $d= \partial_t + \bar\partial_t$ depends on $t$. 

 Throughout the paper, we shall denote $n:=\mbox{dim}_\C X_t$ for $t\in\Delta$.

 The spaces of $C^\infty$ forms of degree $k$, resp. of bidegree $(p,\,q)$, on $X_t$ will be denoted by $C^{\infty}_k(X, \, \C)$, resp. $C^{\infty}_{p, \, q}(X_t, \, \C)$. Given a form $u$, its component of type $(p, q)$ with respect to the complex structure $J_t$ will be denoted by $u_t^{p,\, q}$. 

 The $\lambda$-eigenspace of a given elliptic differential operator $P_t : C^{\infty}_{p, \, q}(X_t, \, \C)\rightarrow C^{\infty}_{p, \, q}(X_t, \, \C)$ will be denoted by $E^{p, \, q}_{P_t}(\lambda)$, where in most cases $P_t$ will be taken to be one of the Laplace-Beltrami operators $\Delta_t', \,\Delta_t''$ associated with a given Hermitian metric on $X_t$. 

 As usual, we shall denote by $h^{p, \, q}(t):=\mbox{dim}_\C H^{p,\,q}(X_t,\,\C)$, resp. $b_k:=\mbox{dim}_\C H^k(X,\,\C)$ the Hodge, resp. Betti numbers of $X_t$. Thanks to the $\partial\bar\partial$ assumption on $X_t$ for every $t\neq 0$, every $h^{p, \, q}(t)$ is constant on $\Delta^{\star}$ after possibly shrinking $\Delta$ about $0$. However, it may happen that $h^{p, \, q}(0)>h^{p, \, q}(t)$ for $t\neq 0$, although this case is {\it a posteriori} ruled out if $X_t$ is projective for every $t\in\Delta^\star$ by Theorem \ref{The:projdef}. We stress that the weaker $\partial\bar\partial$ assumption for $t\neq 0$ in Theorem \ref{The:limitsG} need not rule out, even a posteriori, the jumping at $t=0$ of $h^{p, \, q}(t)$.

\begin{Rem}\label{Rem:restriction}{\rm One of the main difficulties one is faced with in trying to prove Theorems \ref{The:projdef} and \ref{The:limitsG} is the possible jump at $t=0$ of the Hodge numbers $h^{p, \, q}(t)$.

 Specifically, suppose a $C^{\infty}$ family $(\gamma_t)_{t\in\Delta}$ of Hermitian metrics on the fibres $(X_t)_{t\in\Delta}$ has been fixed. We get an associated $C^{\infty}$ family $(\Delta''_t)_{t\in\Delta}$ of Laplace-Beltrami operators acting on the $J_t$-$(p, \, q)$-forms of $X$ for every bidegree $(p,\,q)$. For every $t\in\Delta$ and every bidegree, the operator $\Delta_t'':=\bar\partial_t\bar\partial_t^\star + \bar\partial_t^\star\bar\partial_t$ is elliptic and therefore has a compact resolvent and a discrete spectrum

\begin{equation}\label{eqn:spectrum}0=\lambda_0(t)\leq \lambda_1(t) \leq \dots \leq \lambda_k(t) \leq  \dots\end{equation} 

\vspace{1ex}

\noindent with $\lambda_k(t)\rightarrow +\infty$ as $k\rightarrow +\infty$. Thanks to the Hodge isomorphism 

$$H^{p,\,q}(X_t,\,\C)\simeq\ker(\Delta_t'' : C^\infty_{p,\,q}(X_t,\,\C)\to C^\infty_{p,\,q}(X_t,\,\C)),  \hspace{3ex} t\in\Delta,$$

\noindent the multiplicity of zero as an eigenvalue of $\Delta''_t$ equals $h^{p, q}(t)$. By results of Kodaira and Spencer (see [Kod85, Lemmas 7.5-7.7 and Proof of Theorem 7.2, p. 338-343]), for every small $\varepsilon>0$, the number $m\in\N^{\star}$ of eigenvalues (counted with multiplicities) of $\Delta''_t$ contained in the interval $[0, \, \varepsilon)$ is independent of $t$ if $t\in\Delta$ is sufficiently close to $0$ (say $\delta_{\varepsilon}$-close). If $\varepsilon>0$ has been chosen so small that $0$ is the only eigenvalue of $\Delta''_0$ contained in $[0, \, \varepsilon)$, it follows that $m=h^{p, q}(0)\geq h^{p, q}(t)$ for $t$ sufficiently close to $0$ (the upper-semicontinuity property). Consequently, for $t$ near $0$, $h^{p, q}(0) = h^{p, q}(t)$ if and only if $0$ is the only eigenvalue of $\Delta''_t$ lying in $[0, \, \varepsilon)$. In other words, if $h^{p, \, q}(0)>h^{p, q}(t)$ when $t (\neq 0)$ is near $0$, choosing increasingly small $\varepsilon >0$ gives eigenvalues of $\Delta''_t$

\begin{equation}\label{eqn:smalleigenvalues} 0<\lambda_{k_1}(t)\leq \lambda_{k_2(t)} \leq \dots \leq \lambda_{k_N}(t):=\varepsilon_t < \varepsilon,  \hspace{3ex}  t\in\Delta^{\star},\end{equation}

\noindent that converge to zero (i.e. $\varepsilon_t\rightarrow 0$) when $t\rightarrow 0$, where $N=h^{p, \, q}(0)-h^{p, \, q}(t)$. 

Now, we will have to solve on several occasions throughout the paper, for $C^{\infty}$ (up to $t=0$) families of $\bar\partial_t$-exact forms $(v_t)_{t\in\Delta}$, equations of the shape

$$\bar\partial_t u_t = v_t \hspace{3ex} \mbox{on} \hspace{1ex} X_t, \hspace{3ex} \mbox{for} \hspace{1ex} t\in\Delta\setminus\{0\},$$

\noindent whose minimal $L^2_{\gamma_t}$-norm solution is given by the Neumann formula

$$u_t = \Delta_t^{''-1}\bar\partial_t^\star v_t, \hspace{3ex}  t\in\Delta\setminus\{0\},$$

\noindent that features the Green operator $\Delta_t^{''-1}$ of $\Delta''_t$ (i.e. the inverse of the restriction of $\Delta''_t$ to the orthogonal complement of its kernel). The inverses $1/\lambda_{k_j}(t)$ of the eigenvalues of $\Delta''_t$ are eigenvalues for $\Delta_t^{''-1}$ and $1/\lambda_{k_j}(t)\rightarrow +\infty$ when $t\rightarrow 0$ for every $k_j\in\{k_1, \dots , k_N\}$ (i.e. for every {\it small} eigenvalue) if there is a jump $h^{p,\,q}(0)>h^{p,\,q}(t)$. It follows that, if $\bar{\partial}_t^{\star}v_t$ has non-trivial projections onto the eigenspaces $E^{p, \, q}_{\Delta''_t}(\lambda_{k_j(t)})$ with $k_j\in\{k_1, \dots , k_N\}$, these projections get multiplied by $1/\lambda_{k_j}(t)$ when $\Delta_t^{''-1}$ is applied to $\bar{\partial}_t^{\star}v_t$. Then $u_t$ need not be bounded as $t$ approaches $0$, unless the said projections can be proved to tend to zero sufficiently quickly to offset the growth $1/\lambda_{k_j}(t)\rightarrow +\infty$ when $t$ approaches $0$. This unboundedness may cause the family of forms $(u_t)_{t\in\Delta\setminus\{0\}}$ to not extend across $t=0$, i.e. to not have a limit $u_0$ on $X_0$ when $t\rightarrow 0$. 

 The same conclusion applies to the $\partial_t$-Laplacians $\Delta_t':=\partial_t\partial_t^\star + \partial_t^\star\partial_t$ because of the possible jump (upwards) at $t=0$ of the dimensions of the cohomology groups $H^{p,\,q}_{\partial_t}(X_t,\C)$ that depend on the complex structures $J_t$.

 }

\end{Rem}

 \begin{Rem}\label{Rem:restriction_De-Rham}{\rm However, a simple observation that will play a major role in our approach to Theorem \ref{The:limitsG} is that the unboundedness phenomenon described in Remark \ref{Rem:restriction} does not occur for the $C^{\infty}$ family $(\Delta_t)_{t\in\Delta}$ of $d$-Laplacians on the fibres $(X_t)_{t\in\Delta}$.

Indeed, analogous families $(u_t=\Delta_t^{-1}d^\star_t v_t)_{t\in\Delta\setminus\{0\}}$ of minimal $L^2_{\gamma_t}$-norm solutions of $d$-equations

$$d u_t = v_t \hspace{3ex} \mbox{on} \hspace{1ex} X_t, \hspace{3ex}  t\in\Delta\setminus\{0\},$$

\noindent for given $C^{\infty}$ (up to $t=0$) families of $d$-exact forms $(v_t)_{t\in\Delta}$, always extend smoothly to a form $u_0$ on $X_0$ thanks to the De Rham cohomology of the fibres $X_t$ being independent of $t\in\Delta$. The reason is that the family of manifolds $(X_t)_{t\in\Delta}$ is $C^{\infty}$ trivial, so the Betti numbers $b_k$ of the fibres $X_t$ are constant. Therefore, there is no ``jumping'' phenomenon in this case.  
 
  }

\end{Rem}

\section {Preliminaries to the proof of Theorem \ref{The:limitsG}}\label{section:preliminaries_proof} 

We now start the proof of Theorem \ref{The:limitsG} that will occupy the rest of the paper.

\subsection {Reduction of the uniform boundedness problem to a positivity problem}\label{subsection:reduction_positivity}

Fix any $C^{\infty}$ family $(\gamma_t)_{t\in\Delta}$ of Gauduchon metrics on the respective fibres $(X_t)_{t\in\Delta}$. (It is well known that such families exist, see a version of the argument in [Pop13, $\S.3$].) For every $k\in\{0,\dots ,2n\}$ and every $p,q\in\{0,\dots ,n\}$, we denote by $\Delta_t:= dd_t^\star + d_t^\star d : C^\infty_k(X,\C)\to C^\infty_k(X,\C)$ and by 

$$\Delta_t':=\partial_t\partial_t^\star + \partial_t^\star\partial_t, \,\,\, \Delta_t'':=\bar\partial_t\bar\partial_t^\star + \bar\partial_t^\star\bar\partial_t : C^\infty_{p,\,q}(X_t,\C) \to C^\infty_{p,\,q}(X_t,\C)$$

\noindent the $d$-, $\partial_t$- and $\bar\partial_t$-Laplace-Beltrami operators induced by the metrics $\gamma_t$ on $X_t$. Let $(\lambda_j(t))_{j\in\N}$ be the eigenvalues, ordered non-increasingly and repeated as many times as the respective multiplicity, of 

$$\Delta''_t:C^{\infty}_{n, n-1}(X_t, \C)\longrightarrow C^{\infty}_{n, n-1}(X_t, \C), \hspace{3ex} t\in\Delta.$$ 

\noindent By [Kod86], each $\lambda_j$ is a continuous function of $t\in\Delta$. If there are eigenvalues such that $\lambda_j(t)>0$ for $t\neq 0$ and $\lambda_j(0)= 0$, there are only finitely many of them numbering $h^{n, \, n-1}(0)-h^{n, \, n-1}(t)=h^{0, 1}(0)-h^{0, 1}(t)$ for any $t\neq 0$ close to $0$. This number is, of course, independent of $t\neq 0$. For $t\neq 0$, let $\varepsilon_t''>0$ be the largest of these {\it small} eigenvalues, so $\varepsilon_t''\rightarrow 0$ as $t\rightarrow 0$. The remaining, infinitely many, eigenvalues are then bounded below (after possibly shrinking $\Delta$ about $0$) by some $\varepsilon''>0$ independent of $t\in\Delta$. Thus 

\begin{equation}\label{eqn:spec}\mbox{Spec}\,\Delta''_t\subset [0,\,\, \varepsilon_t''] \cup [\varepsilon'',\,\, +\infty), \hspace{3ex} t\in\Delta,\end{equation}

\noindent where we have set $\varepsilon_0''=0$. We get an $L^2_{\gamma_t}$-orthogonal eigenspace decomposition 

\begin{equation}\label{eqn:orthdecomp}C^{\infty}_{n, n-1}(X_t, \C)=\bigoplus\limits_{\lambda\leq\varepsilon_t''}E^{n, n-1}_{\Delta''_t}(\lambda) \oplus \bigoplus\limits_{\lambda\geq\varepsilon''}E^{n, n-1}_{\Delta''_t}(\lambda), \hspace{3ex} t\in\Delta.\end{equation}

\noindent Now, $\Delta''_t$ being an elliptic self-adjoint operator, it has a compact resolvent and there exists an orthonormal basis $(e_j^{n, n-1}(t))_{j\in\N}$ of $C^{\infty}_{n, n-1}(X_t, \C)$ consisting of eigenvectors of $\Delta''_t$ :

\begin{equation}\label{eqn:eigen}\Delta''_t\, e_j^{n, n-1}(t)= \lambda_j(t)\, e_j^{n, n-1}(t), \hspace{3ex} t\in\Delta.\end{equation} 

\noindent Furthermore, in the three-space orthogonal decomposition 

\begin{equation}\label{eqn:3spacedecomp}C^{\infty}_{n, n-1}(X_t, \C)=\ker\Delta''_t \oplus \mbox{Im}\,\bar\partial_t \oplus \mbox{Im}\,\bar\partial_t^{\star},\end{equation}

\noindent each subspace is $\Delta''_t$-invariant due to $\Delta''_t$ commuting with $\bar\partial_t$ and $\bar\partial^{\star}_t$. This means that the eigenvectors $e_j^{n, n-1}(t)$ forming an orthonormal basis can be chosen such that each of them lies in one (and only one) of the three subspaces of (\ref{eqn:3spacedecomp}). So none of the $e_j^{n, n-1}(t)$ straddles two or three subspaces. These simple reductions are valid for every $t\in\Delta$ and we will henceforth suppose that the choices have been made as described above. The orthogonal decomposition of $\partial_t\gamma_t^{n-1}\in C^{\infty}_{n, n-1}(X_t, \C)$ according to (\ref{eqn:orthdecomp}) has the shape:

\begin{equation}\label{eqn:del-decomp}\partial_t\gamma_t^{n-1}= \sum\limits_{j\in J_1}c_j(t)\, e_j^{n, n-1}(t) + \sum\limits_{j\in J_2}c_j(t)\, e_j^{n, n-1}(t) = U_t + V_t, \hspace{3ex} t\in\Delta,\end{equation}

\noindent where $U_t = \sum\limits_{j\in J_1}c_j(t)\, e_j^{n, n-1}(t)\in\bigoplus\limits_{\lambda\leq\varepsilon_t''}E^{n, n-1}_{\Delta''_t}(\lambda)$ and $V_t = \sum\limits_{j\in J_2}c_j(t)\, e_j^{n, n-1}(t)\in\bigoplus\limits_{\lambda\geq\varepsilon''}E^{n, n-1}_{\Delta''_t}(\lambda)$, with coefficients $c_j(t)\in\C^{\star}$ and index sets $J_1, J_2\subset\N$ such that $J_1\cap J_2=\emptyset$. As already noticed, by the Gauduchon condition, $\partial_t\gamma_t^{n-1}$ is $d$-closed for all $t\in\Delta$ and, since it is $\partial_t$-exact, it must also be $\bar\partial_t$-exact for all $t\neq 0$ by the $\partial\bar\partial$-lemma. Since each eigenvector $e_j^{n, n-1}(t)$ belongs to one of the three orthogonal subspaces of (\ref{eqn:3spacedecomp}), this means that only eigenvectors belonging to $\mbox{Im}\,\bar\partial_t$ can have a non-trivial contribution to (\ref{eqn:del-decomp}) for $t\neq 0$.

 In particular, for every $t\neq 0$, both $U_t$ and $V_t$ are $\bar\partial_t$-exact. We can therefore find, for every $t\neq 0$, a smooth $J_t\!-\!(n, n-2)$-form $w_t$ such that $V_t=\bar\partial w_t$. If we choose the form $w_t$ of minimal $L^2$ norm (with respect to $\gamma_t$) with this property, the condition $V_t\in\bigoplus\limits_{\lambda\geq\varepsilon''}E^{n, n-1}_{\Delta''_t}(\lambda)$ guarantees that the family of forms $(w_t)_{t\in\Delta^{\star}}$ extends smoothly across $t=0$ to a family $(w_t)_{t\in\Delta}$ varying in a $C^{\infty}$ way with $t$ up to $t=0$. This is because the eigenvalues $\lambda$ contributing to $V_t$ are uniformly bounded below by $\varepsilon''>0$.

 As for $U_t\in\bigoplus\limits_{\lambda\leq\varepsilon_t''}E^{n, n-1}_{\Delta''_t}(\lambda)$, we are unable to guarantee the boundedness near $t=0$ of its $\bar\partial_t$-potential because of the eigenvalues $\lambda_j(t)\leq\varepsilon_t''$ converging to $0$. Therefore we will not consider the $\bar\partial_t$-potential. However, the $(n, n-1)$-form $U_t$ is $d$-closed. Indeed, it is $\partial_t$-closed in a trivial way for bidegree reasons and is also $\bar\partial_t$-closed (even $\bar\partial_t$-exact, as it has been argued above). Thus, the $\partial\bar\partial$-lemma implies that $U_t$ is $d$-exact for every $t\neq 0$. We can therefore find, for all $t\neq 0$, a form $\xi_t$ of degree $2n-2$ such that $U_t=d\,\xi_t$. If we choose the form $\xi_t$ of minimal $L^2$-norm (with respect to $\gamma_t$) with this property, we have

\begin{equation}\label{eqn:dmin}\xi_t=\Delta_t^{-1}\, d_t^{\star}U_t, \hspace{3ex} t\neq 0,\end{equation}

\noindent where, for all $t\in\Delta$ (including $t=0$), $\Delta_t=d\, d_t^{\star} + d_t^{\star}\, d : C^{\infty}_{2n-2}(X, \, \C)\rightarrow C^{\infty}_{2n-2}(X, \, \C)$ is the $d$-Laplacian associated with the metric $\gamma_t$ and $\Delta_t^{-1}$ is the inverse of the restriction of $\Delta_t$ to the orthogonal complement of its kernel (the Green operator of $\Delta_t$). Now, the Hodge isomorphism theorem gives:

\begin{equation}\label{eqn:Hodgeiso-d}\ker\Delta_t\simeq H_{DR}^{2n-2}(X_t, \C)=H^{2n-2}(X, \C), \hspace{3ex} t\in\Delta,\end{equation}

\noindent and we know that all the De Rham cohomology groups $H_{DR}^{2n-2}(X_t, \C)$ of the fibres $X_t$ can be identified with a fixed space $H^{2n-2}(X, \C)$. In particular, the dimension of $\ker\Delta_t$ is independent of $t\in\Delta$, which means that the positive eigenvalues of $\Delta_t$ have a uniform positive ($>0$) lower bound for $t$ close to $0$ (cf. Kodaira-Spencer arguments [Kod86] recalled in Remarks \ref{Rem:restriction} and \ref{Rem:restriction_De-Rham} and applied to the $C^{\infty}$ family of strongly elliptic operators $(\Delta_t)_{t\in\Delta}$). Thus, in this respect, there is a sharp contrast between the $d$-Laplacian $\Delta_t$ and its $\bar\partial_t$-counterpart $\Delta''_t$: unlike $\Delta''_t$, $\Delta_t$ never displays the small eigenvalue phenomenon. In particular, the family of $(2n-2)$-forms $(\xi_t)_{t\in\Delta^{\star}}$ extends smoothly across $t=0$ to a family $(\xi_t)_{t\in\Delta}$ of forms varying in a $C^{\infty}$ way with $t\in\Delta$ (up to $t=0$). 

\vspace{1ex}

 Our discussion so far can be summed up as follows.

\begin{Lem}\label{Lem:summingup} Given any family of Gauduchon metrics $(\gamma_t)_{t\in\Delta}$ varying in a $C^{\infty}$ way with $t\in\Delta$ on the fibres of a family $(X_t)_{t\in\Delta}$ in which the $\partial\bar{\partial}$-lemma holds on $X_t$ for every $t\neq 0$, we can find a decomposition:

\begin{equation}\label{eqn:half-conclusion1}\partial_t\gamma_t^{n-1}=d\, \xi_t + \bar\partial_t w_t, \hspace{3ex} t\in\Delta,\end{equation}

\noindent in such a way that

\begin{equation}\label{eqn:half-conclusion2}d\, \xi_t\in\bigoplus\limits_{\lambda\leq\varepsilon_t''}E^{n, n-1}_{\Delta''_t}(\lambda), \hspace{3ex} \bar\partial_t w_t\in\bigoplus\limits_{\lambda\geq\varepsilon''}E^{n, n-1}_{\Delta''_t}(\lambda),\end{equation}

\noindent where $(w_t)_{t\in\Delta}$ and $(\xi_t)_{t\in\Delta}$ are families of $(2n-2)$-forms and respectively $(n, n-2)$-forms varying in a $C^{\infty}$ way with $t\in\Delta$ (up to $t=0$), $\varepsilon''>0$ is independent of $t$, $\varepsilon_t''>0$ for $t\neq 0$ and $\varepsilon_t''$ converges to zero as $t$ approaches $0\in\Delta$ (thus $\varepsilon_0''=0$). Moreover, the following identity holds:

\begin{equation}\label{eqn:half-conclusion3}\partial_t(\gamma_t^{n-1}-\xi_t^{n-1, \, n-1})= \bar\partial_t(\xi_t^{n, \, n-2} + w_t), \hspace{3ex} t\in\Delta.\end{equation}

\noindent As the form $\xi_t^{n-1, \, n-1}$ need not be real, we find it more convenient to write:

\begin{equation}\label{eqn:half-conclusion4}\partial_t(\gamma_t^{n-1}-\xi_t^{n-1, \, n-1} - \overline{\xi_t^{n-1, \, n-1}})= \bar\partial_t(\xi_t^{n, \, n-2} + \overline{\xi_t^{n-2, \, n}} + w_t), \hspace{3ex} t\in\Delta.\end{equation}

\end{Lem}

\vspace{2ex}

 To get (\ref{eqn:half-conclusion3}) from (\ref{eqn:half-conclusion1}), it suffices to write $d\, \xi_t = \partial_t\xi_t + \bar\partial_t\xi_t$ and to remember that $d\, \xi_t = U_t$ is a form of pure $J_t$-type $(n, n-1)$. Hence $d\, \xi_t = \partial_t\xi_t^{n-1, \, n-1} + \bar\partial_t\xi_t^{n, \, n-2}$. The vanishing of the $(n-1, n)$-component of $d\, \xi_t$ amounts to $\bar\partial_t\xi_t^{n-1, \, n-1}+\partial_t\xi_t^{n-2, \, n}=0$, or equivalently by conjugation to $\partial_t(-\overline{\xi_t^{n-1, \, n-1}})=\bar\partial_t\overline{\xi_t^{n-2, \, n}}$. Hence (\ref{eqn:half-conclusion4}) follows from (\ref{eqn:half-conclusion3}).  

 As all the forms involved in (\ref{eqn:half-conclusion4}) vary in a $C^{\infty}$ way with $t\in\Delta$ (up to $t=0$), to finish the proof of Theorem \ref{The:projdef} it clearly suffices to show that 

\begin{equation}\label{eqn:posquestion}\gamma_t^{n-1}-\xi_t^{n-1, \, n-1}-\overline{\xi_t^{n-1, \, n-1}}    >0, \hspace{3ex} \mbox{for all} \hspace{2ex} t\in\Delta.\end{equation}

 Indeed, if this positivity property has been proved, Michelsohn's observation in linear algebra [Mic83, p. 279-280] enables one to extract the $(n-1)^{st}$ root of $\gamma_t^{n-1}-\xi_t^{n-1, \, n-1}-\overline{\xi_t^{n-1, \, n-1}}$ and to find, for all $t\in\Delta$, a unique $J_t-(1, 1)$-form $\rho_t>0$ such that

\begin{equation}\label{root}\gamma_t^{n-1}-\xi_t^{n-1, \, n-1}-\overline{\xi_t^{n-1, \, n-1}}=\rho_t^{n-1}, \hspace{3ex} t\in\Delta.\end{equation}

\noindent By construction, $\rho_t$ defines a {\it strongly Gauduchon} metric on $X_t$ for every $t\in\Delta$ thanks to (\ref{eqn:half-conclusion4}). In particular, $X_0$ is a {\it strongly Gauduchon} manifold and Theorem \ref{The:limitsG} follows. It actually suffices to prove (\ref{eqn:posquestion}) for $t=0$.

 Moreover, it would clearly suffice to prove the stronger property:  

\begin{equation}\label{eqn:zerolim}\xi_0^{n-1, n-1}=0.\end{equation}

\noindent If this has been proved, then identity (\ref{eqn:half-conclusion3}) applied to $t=0$ reads $\partial_0\gamma_0^{n-1} = \bar\partial_0(\xi_0^{n, n-2} + w_0)$, hence $\gamma_0$ is a {\it strongly Gauduchon} metric on $X_0$ and Theorem \ref{The:limitsG} follows. 

\vspace{2ex} 

 We have thus reduced our uniform boundedness problem for the {\it main quantity} $I_t$ to the positivity problem (\ref{eqn:posquestion}) or the vanishing subproblem (\ref{eqn:zerolim}).

\vspace{2ex}

\noindent {\it The positivity problem} \\

 Let $||\cdot||=||\cdot||_t$ and $\langle\langle\,\, , \,\, \rangle\rangle = \langle\langle\,\, , \,\, \rangle\rangle_t$ stand for the $L^2$-norm and respectively the $L^2$-scalar product defined by the Gauduchon metric $\gamma_t$ on the forms of $X_t$.

 For the sake of perspicuity, we begin by proving (\ref{eqn:zerolim}) in a special case that brings out the mechanism and locates the difficulty. Different arguments will subsequently be given to settle the positivity problem in full generality.

\subsection {Proof of (\ref{eqn:zerolim}) and implicitly of Theorem \ref{The:limitsG} in an ideal case}\label{subsection:ideal-case}

 Consider the orthogonal decompositions of $\gamma_t^{n-1}$ analogous to (\ref{eqn:del-decomp}) with respect to the eigenspaces of $\Delta_t''$ and respectively $\Delta_t'$ acting on $J_t$-type $(n-1, n-1)$-forms:

\begin{equation}\label{eqn:Delta''decomp}\gamma_t^{n-1}=u_t + v_t, \hspace{2ex} \mbox{with}\,\, u_t\in\bigoplus\limits_{\mu\leq\delta_t}E^{n-1, n-1}_{\Delta_t''}(\mu), \hspace{2ex} v_t\in\bigoplus\limits_{\mu\geq\delta}E^{n-1, n-1}_{\Delta_t''}(\mu), \hspace{3ex} t\in\Delta,\end{equation}

\noindent and

\begin{equation}\label{eqn:Delta'decomp}\gamma_t^{n-1}=\bar{u_t} + \bar{v_t}, \hspace{2ex} \mbox{with}\,\, \bar{u_t}\in\bigoplus\limits_{\mu\leq\delta_t}E^{n-1, n-1}_{\Delta_t'}(\mu), \hspace{2ex} \bar{v_t}\in\bigoplus\limits_{\mu\geq\delta}E^{n-1, n-1}_{\Delta_t'}(\mu), \hspace{3ex} t\in\Delta,\end{equation}

\noindent where $\mbox{Spec}\, (\Delta_t'':C^{\infty}_{n-1, \, n-1}(X_t, \, \C)\rightarrow C^{\infty}_{n-1, \, n-1}(X_t, \, \C))\subset [0, \, \delta_t] \cup [\delta, \, +\infty)$ and $\delta_t\rightarrow 0$ as $t\rightarrow 0$, while $\delta>0$ is independent of $t$. The rest of the notation is analogous to that used earlier for $(n, n-1)$-forms, the symbol $E^{p, q}(\lambda)$ denoting eigenspaces of $(p, q)$-forms with eigenvalue $\lambda$. Decompositions (\ref{eqn:Delta''decomp}) and (\ref{eqn:Delta'decomp}) are conjugate to each other because $\gamma_t^{n-1}$ is a real form and $\overline{\Delta'_t}=\overline{\partial_t\partial_t^{\star} + \partial_t^{\star}\partial_t}= \bar\partial_t\bar\partial_t^{\star} + \bar\partial_t^{\star}\bar\partial_t=\Delta''_t$. In particular, the eigenvalues of both $\Delta'_t$ and $\Delta''_t$ being real (even non-negative, by self-adjointness), the equivalence holds: $u\in E^{n-1, n-1}_{\Delta''_t}(\lambda) \Leftrightarrow \bar{u}\in E^{n-1, n-1}_{\Delta'_t}(\lambda)$.

\begin{Def}\label{Def:idealcase} We say that the {\bf ideal case} occurs if \\

\noindent (i)\, $u_t=\bar{u_t}$ for all $t\in\Delta.$ In other words, the forms $u_t$ and $v_t$ into which $\gamma_t^{n-1}$ splits in (\ref{eqn:Delta''decomp}) are real;  \\

\noindent (ii)\, $\partial_t\Delta''_t u=\Delta''_t\partial_t u$ for all $u\in C^{\infty}_{n-1, \, n-1}(X_t, \C)$ and all $t\in\Delta.$ In other words, $\partial_t$ commutes with $\Delta''_t$ on $J_t$-type $(n-1, n-1)$-forms.

\end{Def}

\vspace{2ex}

 If the Laplacians $\Delta'_t$ and $\Delta''_t$ were calculated with respect to a {\it K\"ahler metric}, then $\Delta'_t=\Delta''_t$ and the {\it ideal case} would occur since $\partial_t$ always commutes with $\Delta_t'$. The failure of the {\it ideal case} to occur in general is caused by the failure of the Gauduchon metric $\gamma_t$ to be K\"ahler.

\begin{Lem}\label{Lem:missinglinkest} Let $(X_t)_{t\in\Delta}$ be any family such that the $\partial\bar\partial$-lemma holds on $X_t$ for every $t\neq 0$. Let $(\gamma_t)_{t\in\Delta}$ be any family of Gauduchon metrics varying in a $C^{\infty}$ way with $t\in\Delta$ on the fibres $(X_t)_{t\in\Delta}$. Suppose the {\bf ideal case} occurs. The notation being that of Lemma \ref{Lem:summingup}, the following estimate holds:

\begin{equation}\label{eqn:pos-estimate}||\xi_t^{n-1, \, n-1}||\leq\varepsilon_0(t)\, ||\gamma_t^{n-1}||, \hspace{3ex} \mbox{for all}\,\, t\in\Delta\setminus\{0\},\end{equation}

\noindent with a constant $\varepsilon_0(t)>0$ converging to zero as $t\rightarrow 0$. In particular, $\xi_0^{n-1, n-1}=0$ and the metric $\gamma_0$ is strongly Gauduchon, proving Theorem \ref{The:limitsG}.

\end{Lem}

 To infer the second statement from estimate (\ref{eqn:pos-estimate}), it suffices to remember that $(\gamma_t^{n-1})_{t\in\Delta}$ and the norms $(||\cdot || = ||\cdot ||_t)_{t\in\Delta}$ vary in a $C^{\infty}$ way with $t$ (up to $t=0$). This clearly implies that $||\gamma_t^{n-1}||$ has a positive upper bound independent of $t$ if $t$ is close to $0$ (it actually converges to $||\gamma_0^{n-1}||\in (0, \, +\infty)$ when $t\rightarrow 0$). Hence $\xi_0^{n-1, n-1}=0$, i.e. (\ref{eqn:zerolim}). 

\vspace{2ex}

\noindent {\it Proof of Lemma \ref{Lem:missinglinkest}.} Condition $(ii)$ of Definition \ref{Def:idealcase} implies that the decomposition (\ref{eqn:half-conclusion1}) of $\partial_t\gamma_t^{n-1}$ (which is known to satisfy (\ref{eqn:half-conclusion2}) and to be unique with this property) is obtained by applying $\partial_t$ on both sides of decomposition (\ref{eqn:Delta''decomp}) of $\gamma_t^{n-1}$. Thus $\varepsilon_t''=\delta_t$, $\varepsilon''=\delta$ and 

\begin{equation}\label{eqn:ii}\partial_t u_t=d\,\xi_t, \hspace{2ex} t\in\Delta.\end{equation}

\noindent On the other hand, property $(i)$ of Definition \ref{Def:idealcase}, combined with (\ref{eqn:Delta'decomp}), gives:

\begin{equation}\label{eqn:i} u_t=\bar{u_t}\in\bigoplus\limits_{\mu\leq\delta_t=\varepsilon_t''}E^{n-1, n-1}_{\Delta_t'}(\mu), \hspace{2ex} t\in\Delta.\end{equation}

 As explained earlier, the positive eigenvalues of the $d$-Laplacian $\Delta_t : C^{\infty}_{2n-2}(X, \, \C)\rightarrow C^{\infty}_{2n-2}(X, \, \C)$ defined by the metric $\gamma_t$ have a positive lower bound independent of $t$ if $t\in\Delta$ is close to $0$. This means that there exists a constant $c>0$, independent of $t\in\Delta$, such that the restriction of $\Delta_t$ to the orthogonal complement of its kernel satisfies

\begin{equation}\label{eqn:Deltalbound}(\Delta_t)_{|(\ker\Delta_t)^{\perp}}\geq c\, \mbox{Id}, \hspace{3ex} t\in\Delta,\end{equation}

\noindent after possibly shrinking the base $\Delta$ about $0$. Putting the bits together, we get the following estimate:

\begin{eqnarray}\label{eqn:estseq}\nonumber c\,||\xi_t^{n-1, \, n-1}||^2 & \leq & c\,||\xi_t||^2\leq \langle\langle \Delta_t\xi_t, \, \xi_t\rangle\rangle=||d\,\xi_t||^2 \\
\nonumber  & = & ||\partial_t u_t||^2\leq\langle\langle \Delta'_t u_t, \, u_t\rangle\rangle \leq \varepsilon_t''\, ||u_t||^2 \\
  & \leq & \varepsilon_t''\, ||\gamma_t^{n-1}||^2, \hspace{3ex} \mbox{for all}\,\, t\in\Delta\setminus\{0\}.\end{eqnarray}

 \noindent Indeed, on the first line: the first inequality follows from the components $\xi_t^{n, \, n-2}$, $\xi_t^{n-1, \, n-1}$, $\xi_t^{n-2, \, n}$, that split $\xi_t$ into $J_t$-types, being mutually orthogonal as forms of different pure types; the second inequality follows from (\ref{eqn:Deltalbound}) as $\xi_t$ has been chosen of minimal $L^2$-norm in (\ref{eqn:dmin}), hence $\xi_t\in\mbox{Im}\, d_t^{\star}$, so, in particular, $\xi_t\in(\ker\Delta_t)^{\perp}$ and (\ref{eqn:Deltalbound}) applies; the identity follows from $d_t^{\star}\xi_t=0$ which holds because $\mbox{Im}\, d_t^{\star}\subset\ker d_t^{\star}$. Further down on the second line: the equality with the last term of the first line follows from (\ref{eqn:ii}); the first inequality is obvious as $\langle\langle \Delta'_t u_t, \, u_t\rangle\rangle = ||\partial_t u_t||^2 + ||\partial_t^{\star} u_t||^2$; the second inequality follows from (\ref{eqn:i}). Finally, the inequality between the last term of the second line and the term on the third line is obvious from the decomposition (\ref{eqn:Delta''decomp}) being orthogonal. The conclusion of (\ref{eqn:estseq}) is that

\begin{equation}\label{eqn:icconcl}||\xi_t^{n-1, \, n-1}||^2\leq\frac{\varepsilon_t''}{c}\, ||\gamma_t^{n-1}||^2, \hspace{3ex} \mbox{for all}\,\, t\in\Delta\setminus\{0\},\end{equation}

\noindent which is nothing but estimate (\ref{eqn:pos-estimate}) with $\varepsilon_0(t):=\frac{\varepsilon_t''}{c}\rightarrow 0$ as $t\rightarrow 0$ that we had set out to prove. This concludes the proof of Lemma \ref{Lem:missinglinkest}.  \hfill $\Box$

\vspace{1ex}

 Notice that if we disregard estimate (\ref{eqn:pos-estimate}), the weaker conclusion $d\,\xi_0=0$, which suffices for our purposes since it gives $\partial_0\gamma_0^{n-1}=\bar\partial_0 w_0$ hence $\gamma_0$ is {\it strongly Gauduchon}, can be reached by a quicker route. Indeed, by (\ref{eqn:half-conclusion2}), $d\,\xi_t$ lies in the $\Delta_t''$-eigenspaces with eigenvalues $\lambda\leq\varepsilon_t''\rightarrow 0$ as $t\rightarrow 0$. Hence $d\,\xi_0$ is $\Delta_0''$-harmonic (or, equivalently, both $\bar\partial_0$ and $\bar\partial_0^{\star}$-closed). Moreover, if the {\it ideal case} occurs, (\ref{eqn:ii}) and (\ref{eqn:i}) show that $d\,\xi_t$ has the similar property with respect to $\Delta_t'$, hence $d\,\xi_0$ is also $\Delta_0'$-harmonic (or, equivalently, both $\partial_0$ and $\partial_0^{\star}$-closed). Now, since $d\,\xi_0$ is of pure type and harmonic for both $\Delta_0'$ and $\Delta_0''$, it must be $\Delta_0$-harmonic (i.e. both $d$ and $d_0^{\star}$-closed) since $d=\partial_0 + \bar\partial_0$ and $d_0^{\star}=\partial_0^{\star} + \bar\partial_0^{\star}$. As $d\,\xi_0$ is obviously $d$-exact and as the spaces $\ker\Delta_0$ and $\mbox{Im}\, d$ are orthogonal, we must have $d\,\xi_0=0$.

\subsection{Sufficiency of a small $L^2$ norm for the correcting form}\label{subsection:L2-norm_sufficiency}

 As $\gamma_t^{n-1}>0$, we shall now see that in order to prove Theorem \ref{The:limitsG}, it suffices to show that the $L^2$-norm $||\cdot ||$ of $\xi_t^{n-1, \, n-1}$ can be made arbitrarily small (hence so can the $L^2$-norm of the real form $\xi_t^{n-1, \, n-1} + \overline{\xi_t^{n-1, \, n-1}}$) uniformly w.r.t. $t\in\Delta$. It would actually suffice to guarantee this property when $t=0$ as the following observation shows.

\begin{Lem}\label{Lem:small-correction} Suppose that for an $\varepsilon>0$ independent of $t\in\Delta$, we have

\begin{equation}\label{eqn:small-correction}||\xi_t^{n-1, \, n-1}||<\varepsilon \hspace{3ex} \forall t\in\Delta, \hspace{3ex} \mbox{or merely} \hspace{2ex} ||\xi_0^{n-1, \, n-1}||<\varepsilon.\end{equation}

 Then, if $\varepsilon$ is sufficiently small, there exists a $C^{\infty}$ form $\rho_0 >0$ that is positive definite and of type $(1, \, 1)$ for $J_0$ such that

\begin{equation}\label{eqn:correction-final}\partial_0\rho_0^{n-1} - \partial_0\bigg(\gamma_0^{n-1}-\xi_0^{n-1, \, n-1} - \overline{\xi_0^{n-1, \, n-1}}\bigg) \in \mbox{Im}\, (\partial_0\bar\partial_0).\end{equation}

\noindent In particular, since $\partial_0(\gamma_0^{n-1}-\xi_0^{n-1, \, n-1} - \overline{\xi_0^{n-1, \, n-1}})$ is $\bar\partial_0$-exact by (\ref{eqn:half-conclusion4}), we see that $\partial_0\rho_0^{n-1}$ is $\bar\partial_0$-exact, hence $\rho_0$ is a {\bf strongly Gauduchon} metric on $X_0$.

\end{Lem}

\noindent {\it Proof.} To lighten the notation, we drop the indices and spell out the argument on an arbitrary compact complex $n-$fold $X$ which will be taken to be $X_0$ in the end. Recall that for all $(p, \, q)$, the Aeppli cohomology groups, which are dual to the Bott-Chern cohomology groups $H^{n-p, \, n-q}_{BC}(X, \C)$, are defined by

$$H^{p, \, q}_A(X, \, \C)=\frac{\ker(\partial\bar\partial:C^{\infty}_{p, \, q}(X)\rightarrow C^{\infty}_{p+1, \, q+1}(X))}{\mbox{Im}(\partial:C^{\infty}_{p-1, \, q}(X)\rightarrow C^{\infty}_{p, \, q}(X)) + \mbox{Im}(\bar\partial:C^{\infty}_{p, \, q-1}(X)\rightarrow C^{\infty}_{p, \, q}(X))}.$$

 We refer to the work [Sch07, 2.c., p. 9-10] of M. Schweitzer for the following set-up. Having fixed the metric $\gamma (=\gamma_0)$ on $X (=X_0)$ and calculating all the formal adjoint operators w.r.t. $\gamma$, the fourth-order Laplacian (cf. [KS60])

$$\Delta^{p, \, q}_A:=\partial\partial^{\star} + \bar\partial\bar\partial^{\star} + (\partial\bar\partial)^{\star}(\partial\bar\partial) : C^{\infty}_{p, \, q}(X, \, \C)\longrightarrow C^{\infty}_{p, \, q}(X, \, \C),$$

\noindent which is quite natural to consider in relation to the Aeppli cohomology, is not elliptic, but it can be modified to an elliptic fourth-order Laplacian $\tilde\Delta^{p, \, q}_A : C^{\infty}_{p, \, q}(X, \, \C)\longrightarrow C^{\infty}_{p, \, q}(X, \, \C)$ having the same kernel:

$$\tilde\Delta^{p, \, q}_A:=\partial\partial^{\star} + \bar\partial\bar\partial^{\star} + (\partial\bar\partial)^{\star}(\partial\bar\partial) + (\partial\bar\partial)(\partial\bar\partial)^{\star} + \partial\bar\partial^{\star}(\partial\bar\partial^{\star})^{\star} + (\partial\bar\partial^{\star})^{\star}\partial\bar\partial^{\star}.$$

 By ellipticity, we get the following three-space decomposition which is orthogonal w.r.t. the $L^2$ inner product defined by $\gamma (=\gamma_0)$ on $X (=X_0)$:

$$C^{\infty}_{p, \, q}(X, \C)=\ker\tilde\Delta^{p, \, q}_A \oplus (\mbox{Im}\partial + \mbox{Im}\bar\partial) \oplus \mbox{Im}(\partial\bar\partial)^{\star},$$

\noindent the orthogonal direct sum of the first two subspaces being the kernel of $\partial\bar\partial$:

\begin{equation}\label{eqn:A-kernel}\ker(\partial\bar\partial)=\ker\tilde\Delta^{p, \, q}_A \oplus (\mbox{Im}\partial + \mbox{Im}\bar\partial),\end{equation}

\noindent a decomposition proving the Hodge isomorphism $H^{p, \, q}_A(X, \, \C)\simeq \ker\tilde\Delta^{p, \, q}_A$.

 Taking $(p, \, q)=(n-1, \, n-1)$ in this general context described in [Sch07], recall that we have (cf. (\ref{eqn:half-conclusion3}) at $t=0$ with indices dropped, set $\xi:=\xi_0$):

$$\partial(\gamma^{n-1}-\xi^{n-1, \, n-1}) = \bar\partial(\xi^{n, \, n-2}+w).$$

\noindent Since $\partial\bar\partial\gamma^{n-1}=0$, taking $\bar\partial$ on both sides of the above identity, we get $\partial\bar\partial\xi^{n-1, \, n-1}=0$, hence the following decomposition according to (\ref{eqn:A-kernel}):

\begin{equation}\label{eqn:xi-decomp}\ker(\partial\bar\partial)\ni\xi^{n-1, \, n-1}=\xi^{n-1, \, n-1}_{\tilde\Delta_A} + (\partial\zeta + \bar\partial\eta),\end{equation}

\noindent where $\zeta$ and $\eta$ are $C^{\infty}$ forms of respective types $(n-2, \, n-1)$ and $(n-1, \, n-2)$, while $\xi^{n-1, \, n-1}_{\tilde\Delta_A}\in\ker\tilde\Delta^{n-1, \, n-1}_A$ is orthogonal onto the sum $\partial\zeta + \bar\partial\eta$. By orthogonality, we get:

\begin{equation}\label{eqn:A-smallness}0\leq ||\xi^{n-1, \, n-1}_{\tilde\Delta_A}||\leq ||\xi^{n-1, \, n-1}||<\varepsilon,\end{equation}

\noindent the last inequality being the hypothesis (\ref{eqn:small-correction}) (for $\xi^{n-1, \, n-1}:=\xi^{n-1, \, n-1}_0$). 

 Thus the $\tilde\Delta^{n-1, \, n-1}_A$-harmonic form $\xi^{n-1, \, n-1}_{\tilde\Delta_A}$ is small in $L^2$-norm by (\ref{eqn:A-smallness}). However, the harmonicity w.r.t. an elliptic operator implies that $\xi^{n-1, \, n-1}_{\tilde\Delta_A}$ must be small in a much stronger norm. Indeed, applying the {\it fundamental a priori inequality} satisfied by elliptic operators to the fourth-order elliptic operator $\tilde\Delta^{n-1, \, n-1}_A$, we get for every $k\in\N$ and every $L^2$-form $u$ of type $(n-1, \, n-1)$ such that $\tilde\Delta^{n-1, \, n-1}_Au$ is in the Sobolev space $W^k(X, \, \Lambda^{n-1, \, n-1}T^{\star}X)$ of $(n-1, \, n-1)$-forms on $X$ whose derivatives up to order $k$ are in $L^2:$

\begin{equation}\label{eqn:apriori-ineq}||u||_{W^{k+4}}\leq C_k(||\tilde\Delta_A^{n-1, \, n-1}u||_{W^k} + ||u||_{L^2}),\end{equation}

\noindent where $||\cdot||_{L^2}:= ||\cdot||$ and $C_k>0$ is a constant depending only on $k$. If $u=\xi^{n-1, \, n-1}_{\tilde\Delta_A}$, $\tilde\Delta^{n-1, \, n-1}_Au=0$ and, by (\ref{eqn:A-smallness}), $||u||_{L^2}<\varepsilon$. Thus (\ref{eqn:apriori-ineq}) reduces to \begin{equation}\label{eqn:S-estimate}||\xi^{n-1, \, n-1}_{\tilde\Delta_A}||_{W^{k+4}}\leq C_k\,\varepsilon, \hspace{3ex} k\in\N.\end{equation}

\noindent Now by the well-known Sobolev Lemma, we have a continuous injection:

$$W^k(X, \Lambda^{n-1, \, n-1}T^{\star}X)\hookrightarrow C^l(X, \Lambda^{n-1, \, n-1}T^{\star}X),  \hspace{3ex} \forall k> l+n,$$

\noindent into the space of $(n-1, \, n-1)$-forms of class $C^l$ on $X$. Choosing $l=0$ and $k+4>n$, we get, for a constant $C'_{k+4}>0$ depending only on $k$:

\begin{equation}\label{eqn:C0-est}||\xi^{n-1, \, n-1}_{\tilde\Delta_A}||_{C^0}\leq C'_{k+4}\, ||\xi^{n-1, \, n-1}_{\tilde\Delta_A}||_{W^{k+4}}<C'_{k+4}\, C_k\, \varepsilon,\end{equation}

\noindent having used (\ref{eqn:S-estimate}) for the last inequality. 

 Thus the $C^0$-norm of $\xi^{n-1, \, n-1}_{\tilde\Delta_A}$ can be made arbitrarily small by choosing $\varepsilon$ small enough. Hence so can the $C^0$-norm of $\xi^{n-1, \, n-1}_{\tilde\Delta_A} + \overline{\xi^{n-1, \, n-1}_{\tilde\Delta_A}}$. Since $\gamma^{n-1}>0$, it follows that $\gamma^{n-1} - \xi^{n-1, \, n-1}_{\tilde\Delta_A} - \overline{\xi^{n-1, \, n-1}_{\tilde\Delta_A}}>0$ if $\varepsilon>0$ is chosen small enough, achieving thus the desired positivity property (\ref{eqn:posquestion}). Extracting Michelsohn's $(n-1)^{st}$ root, we get a unique $C^{\infty}$ $(1, 1)$-form $\rho>0$ (i.e. a Hermitian metric $\rho$ on $X$) satisfying

$$\rho^{n-1}=\gamma^{n-1} - \xi^{n-1, \, n-1}_{\tilde\Delta_A} - \overline{\xi^{n-1, \, n-1}_{\tilde\Delta_A}}>0.$$

 On the other hand, it follows from (\ref{eqn:xi-decomp}) that

\begin{eqnarray}\nonumber\gamma^{n-1} - \xi^{n-1, \, n-1} - \overline{\xi^{n-1, \, n-1}} & = & \gamma^{n-1} - \xi^{n-1, \, n-1}_{\tilde\Delta_A} - \overline{\xi^{n-1, \, n-1}_{\tilde\Delta_A}}\\
\nonumber & - & \partial\zeta -\bar\partial\eta -\bar\partial\bar\zeta - \partial\bar\eta \\
\nonumber & = & \rho^{n-1}-\partial\zeta -\bar\partial\eta -\bar\partial\bar\zeta - \partial\bar\eta,\end{eqnarray}

\noindent hence, taking $\partial$ on either side of the above identity, we get

 $$\partial\rho^{n-1} - \partial (\gamma^{n-1} - \xi^{n-1, \, n-1} - \overline{\xi^{n-1, \, n-1}}) = \partial\bar\partial (\eta + \bar\zeta),$$

\noindent proving contention (\ref{eqn:correction-final}) (indices have been dropped here). The proof is complete.  \hfill $\Box$

\section{The iterative procedure and $L^2$ estimates }\label{section:iteration-procedure_L2}

 With Lemma \ref{Lem:small-correction} understood, the rest of the proof of Theorem \ref{The:limitsG} will focus on correcting the forms $\gamma_t^{n-1}>0$ by subtracting real forms whose $L^2$-norms can be made arbitrarily small uniformly w.r.t. $t\in\Delta$ such that the $\partial_t$ of the difference is $\bar\partial_t$-exact for all $t\in\Delta$ ($t=0$ will suffice). We stress that {\it smallness} of the correcting forms in the (relatively weak) $L^2$-norm is enough to achieve the positivity posited in (\ref{eqn:posquestion}) thanks to Lemma \ref{Lem:small-correction}. 

 However, we can see no reason that the $L^2$-norm of $\xi_t^{n-1, \, n-1}$ should be as small as needed if the {\it ideal case} does not occur. In other words, the forms $\xi_t^{n-1, \, n-1}$ constructed in Lemma \ref{Lem:summingup} need not satisfy the hypothesis of Lemma \ref{Lem:small-correction}. Therefore we will replace them by new forms $\widetilde\xi_{t, \, (p)}^{n-1, \, n-1}$ constructed by an inductive procedure that will be described below. The lower index $(p)$ will indicate that $\widetilde\xi_{t, \, (p)}^{n-1, \, n-1}$ has been produced at step $p\in\N$ of the inductive procedure. This procedure is based on an iterative use of Lemma \ref{Lem:summingup} in which $\partial_t\xi_t^{n-1, \, n-1}$ will be replaced by an appropriate form changing at each step $p$. Running the inductive procedure sufficiently many times $p\gg 1$, we shall get the $L^2$-norm $||\widetilde\xi_{t, \, (p)}^{n-1, \, n-1}||$ to become arbitrarily small in a way that is uniform w.r.t. both $t\in\Delta$ and the number $p\gg 1$ of iterations. Uniformity is of the essence in all that follows.

\subsection{The preliminary inductive construction}\label{subsection:preliminary-induction} Nevertheless, an intermediate step is needed in passing from $(\xi_t^{n-1, \, n-1})_{t\in\Delta}$ to $(\widetilde\xi_{t, \, (p)}^{n-1, \, n-1})_{t\in\Delta}$. It consists in running an inductive construction of smooth families of forms $(\xi^{n-1, \, n-1}_{t, \, (p)})_{t\in\Delta}$, $p\in\N$, by iterating Lemma \ref{Lem:summingup} indefinitely. The main observation here is that the $\partial\bar{\partial}$-lemma enables the construction in Lemma \ref{Lem:summingup} to run indefinitely. Identities (\ref{eqn:stepp}) below compare to (\ref{eqn:half-conclusion3}) and (\ref{eqn:stepp'}) to (\ref{eqn:half-conclusion4}).

\begin{Lem}\label{Lem:p-iterations} For every $p\in\N$, there exist families $(\Omega_{t, \, (p)}^{n-1, \, n-1})_{t\in\Delta}$ of $J_t\!-(n-1, \, n-1)-\!$forms varying continuously with $t$, respectively $(\xi_{t, \, (p)})_{t\in\Delta}$ of $(2n-2)$-forms varying in a $C^{\infty}$ way with $t$ (up to $t=0$) such that, for all $t\in\Delta$, we have:

\begin{eqnarray}\label{eqn:stepp}\partial_t(\gamma_t^{n-1}-\Omega_{t, \, (p)}^{n-1, \, n-1}) & = & \partial_t(\gamma_t^{n-1}-\xi_{t, \, (p)}^{n-1, \, n-1}) \\
\nonumber & = & \bar\partial_t(\xi_{t, \, (p)}^{n, \, n-2} + \xi_{t, \, (p-1)}^{n, \, n-2} + \dots + \xi_{t, \, (1)}^{n, \, n-2} + \xi_t^{n, \, n-2} + w_t),\end{eqnarray}

\noindent where, as usual, $\xi_{t, \, (l)}^{r, \, s}$ denotes the component of $J_t$-type $(r, \, s)$ of $\xi_{t, \, (l)}$. As the form $\xi_{t, \, (p)}^{n-1, \, n-1}$ need not be real, we find it more convenient to write:

\begin{equation}\label{eqn:stepp'}\partial_t(\gamma_t^{n-1}-\xi_{t, \, (p)}^{n-1, \, n-1} - \overline{\xi_{t, \, (p)}^{n-1, \, n-1}}) = \bar\partial_t(\xi_{t, \, (p)}^{n, \, n-2} + \overline{\xi_{t, \, (p)}^{n-2, \, n}} + \xi_{t, \, (p-1)}^{n, \, n-2} + \dots + \xi_t^{n, \, n-2} + w_t).\end{equation}

\end{Lem}

\noindent {\it Proof.} Set $\xi_{t, \, (0)}:=\xi_t$ for all $t\in\Delta$. We have already noticed that $\partial_t\gamma_t^{n-1}$ and its projections $d\, \xi_t$ and $\bar\partial_t w_t$ given in (\ref{eqn:half-conclusion1}) are all $d$, $\partial_t$ and $\bar{\partial}_t$-exact for all $t\neq 0$. Writing $d\, \xi_t=\partial_t\xi_t^{n-1, \, n-1} + \bar{\partial}_t\xi_t^{n, \, n-2}$, we see that $\bar\partial_t\xi_t^{n, \, n-2}$ is $\bar{\partial}_t$-closed (even $\bar{\partial}_t$-exact) and is also $\partial_t$-closed for bidegree reasons (being of pure type $(n, n-1)$). Thus $\bar\partial_t\xi_t^{n, \, n-2}$ is $d$-closed and of pure type. By the $\partial\bar{\partial}$-lemma, the $\bar{\partial}_t$-exactness of $\bar\partial_t\xi_t^{n, \, n-2}$ implies its $d$ and $\partial_t$-exactness for all $t\neq 0$. Then $\partial_t\xi_t^{n-1, \, n-1}$ must also be $d$ and $\partial_t$-exact for all $t\neq 0$ as a difference of two such forms. We can thus write

\begin{equation}\label{eqn:iterative-pot}\partial_t\xi_t^{n-1, \, n-1}=\partial_t\Omega_t^{n-1,\, n-1} = d\, \xi_{t,\, (1)}, \hspace{2ex} t\in\Delta,\end{equation}

\noindent where $\Omega_t^{n-1,\, n-1}$ stands for the $\partial_t$-potential of minimal $L^2$-norm $||\cdot ||$ and $\xi_{t, \, (1)}$ denotes the $d$-potential of minimal $L^2$-norm $||\cdot ||$ of $\partial_t\xi_t^{n-1, \, n-1}$. Identities (\ref{eqn:iterative-pot}) {\it a priori} hold only for $t\neq 0$ as the $\partial\bar{\partial}$-lemma is only known to apply on $X_t$ with $t\neq 0$. However, we have seen that in the formula for the minimal $L^2$-norm solution:

\begin{equation}\label{eqn:re-min-form}\xi_{t, \, (1)}=\Delta^{-1}_td^{\star}_t(\partial_t\xi_t^{n-1, \, n-1}), \hspace{3ex} t\in\Delta^{\star},\end{equation}

\noindent the family of Green's operators $(\Delta_t^{-1})_{t\in\Delta}$ is a $C^{\infty}$ family (up to $t=0$) by results of Kodaira-Spencer and the De Rham cohomology being constant on the fibres $X_t$, $t\in\Delta$ (no {\it small eigenvalue} phenomenon for $\Delta_t$). Thus $\partial_0\xi_0^{n-1, \, n-1}$ is $d$-exact and the family $(\xi_{t, \, (1)})_{t\in\Delta}$ is defined and $C^{\infty}$ up to $t=0$.

 Meanwhile, $||\Omega_t^{n-1,\, n-1}||\leq ||\xi_t^{n-1, \, n-1}||$ for all $t\in\Delta^{\star}$ by the $L^2$-norm minimality of $\Omega_t^{n-1,\, n-1}$. As $\xi_t^{n-1, \, n-1}$ is known to extend in a $C^{\infty}$ way to $X_0$, the family $(\Omega_t^{n-1,\, n-1})_{t\in\Delta^{\star}}$ is bounded near $t=0$ and extends at least continuously across $0\in\Delta$ (the eigenvalues of $\Delta'_t$ vary continuously with $t\in\Delta$ by Kodaira-Spencer, see e.g. [Kod85, Theorem 7.2]). Thus identities (\ref{eqn:iterative-pot}) hold for all $t\in\Delta$ (including $t=0$), while the families $(\Omega_t^{n-1,\, n-1})_{t\in\Delta}$ and $(\xi_{t, \, (1)})_{t\in\Delta}$ vary in a continuous, respectively $C^{\infty}$, way with $t$. Set $\Omega_{t, \, 0}^{n-1,\, n-1}\!\!:=\Omega_t^{n-1,\, n-1}$. 

 In view of (\ref{eqn:iterative-pot}), identity (\ref{eqn:half-conclusion3}) becomes:

\begin{equation}\label{eqn:step0}\partial_t(\gamma_t^{n-1}-\Omega_t^{n-1, \, n-1})= \partial_t(\gamma_t^{n-1}-\xi_t^{n-1, \, n-1}) = \bar\partial_t(\xi_t^{n, \, n-2} + w_t), \hspace{2ex} t\in\Delta.\end{equation}

\noindent Writing $d\, \xi_{t, \, (1)}=\partial_t\xi_{t, \, (1)}^{n-1, \, n-1} + \bar\partial_t\xi_{t, \, (1)}^{n, \, n-2}$ (recall that $d\, \xi_{t, \, (1)}$ is of $J_t$-type $(n, \, n-1)$) and using (\ref{eqn:iterative-pot}), we get:

\begin{equation}\label{eqn:step0'}\partial_t(\gamma_t^{n-1}-\xi_{t, \, (1)}^{n-1, \, n-1})=\bar\partial_t(\xi_{t, \, (1)}^{n, \, n-2} + \xi_t^{n, \, n-2} + w_t), \hspace{2ex} t\in\Delta.\end{equation}

 The procedure described above can now be iterated indefinitely. The right-hand term in (\ref{eqn:step0'}) is a $d$-closed and $\bar\partial_t$-exact $(n, n-1)$-form, hence it must be $d$, $\partial_t$ and $\bar{\partial}_t$-exact for all $t\neq 0$ by the $\partial\bar\partial$-lemma. Then so is $\partial_t\xi_{t, \, (1)}^{n-1, \, n-1}$ as a difference of two such forms (i.e. $\partial_t\gamma_t^{n-1}$ and the right-hand term in (\ref{eqn:step0'})). We then get identities analogous to (\ref{eqn:iterative-pot}):

 $$\partial_t\xi_{t, \, (1)}^{n-1, \, n-1} = \partial_t\Omega_{t, \, (1)}^{n-1,\, n-1} = d\, \xi_{t,\, (2)}, \hspace{2ex} t\in\Delta,$$ 

\noindent where $\Omega_{t, \, (1)}^{n-1, \, n-1}$ and $\xi_{t, \, (2)}$ are the $\partial_t$ and respectively $d$-potentials of $\partial_t\xi_{t, \, (1)}^{n-1, \, n-1}$ with minimal $L^2$-norms. They extend continuously, resp. smoothly to $X_0$ by the same arguments as above and, writing $d\, \xi_{t,\, (2)} = \partial_t\xi_{t, \, (2)}^{n-1, \, n-1} + \bar\partial_t\xi_{t, \, (2)}^{n, \, n-2}$, (\ref{eqn:step0'}) reads:

\begin{equation}\label{eqn:step0''}\partial_t(\gamma_t^{n-1}-\xi_{t, \, (2)}^{n-1, \, n-1})=\bar\partial_t(\xi_{t, \, (2)}^{n, \, n-2} + \xi_{t, \, (1)}^{n, \, n-2} + \xi_t^{n, \, n-2} + w_t), \hspace{2ex} t\in\Delta.\end{equation}

 The $(n, \, n-1)$-form $\partial_t\xi_{t, \, (2)}^{n-1, \, n-1}$ is again $d$, $\partial_t$ and $\bar{\partial}_t$-exact for all $t\neq 0$ by the $\partial\bar\partial$-lemma and the procedure can be repeated. At step $p$ one gets:

\begin{equation}\label{eqn:iterative-pot-stepp}\partial_t\xi_{t, \, (p)}^{n-1, \, n-1} = \partial_t\Omega_{t, \, (p)}^{n-1,\, n-1} = d\, \xi_{t,\, (p+1)}, \hspace{2ex} t\in\Delta,\,\, p\in\N,\end{equation} 

\noindent with $\Omega_{t, \, (p)}^{n-1,\, n-1}$ and $\xi_{t,\, (p+1)}$ the $\partial_t$ and respectively $d$-potentials of minimal $L^2$-norms of $\partial_t\xi_{t, \, (p)}^{n-1, \, n-1}$. The form $\Omega_{t, \, (p)}^{n-1,\, n-1}$ can be seen as a correction of $\xi_{t, \, (p)}^{n-1, \, n-1}$ if the latter does not have minimal $L^2$-norm. It is clear that the analogue for $p$ of (\ref{eqn:step0}), (\ref{eqn:step0'}), (\ref{eqn:step0''}) and the definition of $\Omega_{t, \, (p)}^{n-1,\, n-1}$ in (\ref{eqn:iterative-pot-stepp}) add up to the identities (\ref{eqn:stepp}) claimed in the statement. To get (\ref{eqn:stepp'}) from (\ref{eqn:stepp}), recall that $\partial_t\xi_{t, \, (p-1)}^{n-1, \, n-1} = d\, \xi_{t,\, (p)}$ is of $J_t$-type $(n, \, n-1)$, hence its $(n-1, \, n)$-component $\partial_t\xi_{t, \, (p)}^{n-2, \, n} + \bar\partial_t\xi_{t, \, (p)}^{n-1, \, n-1}$ vanishes. Taking conjugates, one gets $\partial_t(-\overline{\xi_{t, \, (p)}^{n-1, \, n-1}})=\bar\partial_t\overline{\xi_{t, \, (p)}^{n-2, \, n}}$ and this term can be added to (\ref{eqn:stepp}) to get (\ref{eqn:stepp'}).   \hfill $\Box$

\vspace{3ex}

 The next observation is that the $L^2$-norm of $\xi_{t, \, (p)}^{n-1, \, n-1}$ can only decrease or stay constant when $p$ increases, so successive iterations of the construction described in Lemma \ref{Lem:p-iterations} bring us increasingly close to achieving our aim of rendering the $L^2$-norm of $\xi_{t, \, (p)}^{n-1, \, n-1}$ arbitrarily small when $p\gg 1$.

\begin{Lem}\label{Lem:decrease} The $C^{\infty}$ families of forms $(\xi_{t, \, (p)}^{n-1, \, n-1})_{t\in\Delta}$, $p\in\N$, constructed in Lemma \ref{Lem:p-iterations} obey the following $L^2$-norm inequalities:

\begin{equation}\label{eqn:decrease}||\xi_{t, \, (p+1)}^{n-1, \, n-1}|| \leq ||\xi_{t, \, (p)}^{n-1, \, n-1}|| \hspace{1ex} \mbox{and} \hspace{1ex} ||\xi_{t, \, (p+1)}|| \leq ||\xi_{t, \, (p)}||    \hspace{3ex}   t\in\Delta, \,\, p\in\N.\end{equation}

\end{Lem}

\noindent {\it Proof.} The minimal $L^2$-norm solutions of equations (\ref{eqn:iterative-pot-stepp}) are given by:

\begin{equation}\label{eqn:re-min-form-p}\xi_{t, \, (p+1)}=\Delta_t^{-1} d_t^{\star} (\partial_t\xi_{t, \, (p)}^{n-1, \, n-1}),  \hspace{3ex}  \mbox{resp.} \hspace{2ex} \Omega_{t, \, (p)}^{n-1, \, n-1}=\Delta_t^{'-1} \partial_t^{\star} (\partial_t\xi_{t, \, (p)}^{n-1, \, n-1}).\end{equation}

 Now it is easily seen that, for any $\partial_t$-exact $(r, \, s)$-form $u$ on $X_t$, one has

\begin{equation}\label{eqn:delta'-halfnorm}||\Delta_t^{'-1}\partial_t^{\star}u||=||\Delta_t^{'-\frac{1}{2}}u||.\end{equation}

\noindent Indeed, if $(e_j^{r, \, s})_{j\in\N}$ is an orthonormal basis of $C^{\infty}_{r, \, s}(X_t, \, \C)$ consisting of eigenvectors of $\Delta'_t$ such that $\Delta'_t e_j^{r, \, s} = \lambda_j\, e_j^{r, \, s}$ and if $u$ splits as $u=\sum\limits_{j\in J_u} c_j \, e_j^{r, \, s}$ with $c_j\in\C$, then $e_j^{r, \, s}$ is $\partial_t$-exact for every $j\in J_u$ and 

\begin{equation}\label{eqn:half-norm-min-sol}\Delta_t^{'-1}\partial_t^{\star} u=\sum\limits_{j\in J_u} \frac{c_j}{\sqrt{\lambda_j}}\, e_j^{r-1, \, s},\end{equation} 

\noindent where $(e_j^{r-1, \, s})_{j\in J_u}$ is an orthonormal subset of $C^{\infty}_{r-1, \, s}(X_t, \, \C)$ consisting of eigenvectors of $\Delta'_t$ corresponding to the same eigenvalues as for $(r, \, s)-$forms: $\Delta'_t e_j^{r-1, \, s} = \lambda_j\, e_j^{r-1, \, s}.$ This is because

$$\partial^{\star}\,\, : \,\, \mbox{Im}\, (\partial :  C^{\infty}_{r-1, \, s}\rightarrow C^{\infty}_{r, \, s}) \longrightarrow \mbox{Im}\, (\partial^{\star} :  C^{\infty}_{r, \, s}\rightarrow C^{\infty}_{r-1, \, s})$$

\noindent is an angle-preserving isomorphism that maps any $\partial$-exact $\Delta'$-eigenvector of type $(r, \, s)$ to a $\Delta'$-eigenvector of type $(r-1, \, s)$ having the same eigenvalue $\lambda$ and an $L^2$-norm multiplied by $\sqrt{\lambda}$. (We have suppressed indices $t$ to ease the notation). A further application of $\Delta^{'-1}$ introduces divisions by the eigenvalues $\lambda_j$, hence the overall effect of applying $\Delta^{'-1}\partial^{\star}$ to $u$ consists in multiplying the coefficients $c_j$ by $\sqrt{\lambda_j}/\lambda_j=1/\sqrt{\lambda_j}$ and replacing the orthonormal set of $(r, \, s)$-forms $\{e_j^{r, \, s}, \, j\in J_u \}$ with an orthonormal set of $(r-1, \, s)$-forms $\{e_j^{r-1, \, s}, \, j\in J_u \}$. Hence (\ref{eqn:half-norm-min-sol}) follows.

On the other hand, $\Delta^{'-\frac{1}{2}}u = \sum\limits_{j\in J_u}\frac{c_j}{\sqrt{\lambda_j}}\, e_j^{r, \, s}$. Thus we get (\ref{eqn:delta'-halfnorm}) since

$$||\Delta_t^{'-1}\partial_t^{\star}u||^2=||\Delta^{'-\frac{1}{2}}u||^2=\sum\limits_{j\in J_u}\frac{|c_j|^2}{\lambda_j}.$$

 Similarly, for any $d$-exact $k$-form $u$ on $X_t$, one has

\begin{equation}\label{eqn:delta-halfnorm}||\Delta_t^{-1}d_t^{\star}u||=||\Delta^{-\frac{1}{2}}u||.\end{equation}

\noindent Thus in the light of (\ref{eqn:re-min-form-p}), (\ref{eqn:delta'-halfnorm}) and (\ref{eqn:delta-halfnorm}) with $u=\partial_t\xi_{t, \, (p)}^{n-1, \, n-1}$, we get

\begin{equation}\label{eqn:half-lap-norms}||\xi_{t, \, (p+1)}||=||\Delta_t^{-\frac{1}{2}}(\partial_t\xi_{t, \, (p)}^{n-1, \, n-1})||, \hspace{3ex} \mbox{resp.} \hspace{2ex} ||\Omega_{t, \, (p)}^{n-1, \, n-1}||=||\Delta_t^{'-\frac{1}{2}}(\partial_t\xi_{t, \, (p)}^{n-1, \, n-1})||.\end{equation}

 We are thus led to compare the Laplacians $\Delta_t'$ and $\Delta_t$ for $t\in\Delta$. We begin by noticing that for any {\it pure-type} (say $(r, \, s)$) form $u$ on some $X_t$, we have:

\begin{equation}\label{eqn:laplace-comparison}\langle\langle\Delta_tu, \, u\rangle\rangle \geq \langle\langle\Delta_t'u, \, u\rangle\rangle,\end{equation}

\noindent while, if $u$ is not $\Delta_t''$-harmonic, we even have

\begin{equation}\label{eqn:laplace-comparison-strict}\langle\langle\Delta_tu, \, u\rangle\rangle > \langle\langle\Delta_t'u, \, u\rangle\rangle.\end{equation}

 Indeed, by compactness of $X_t$, any $(r, \, s)$-form $u$ satisfies: \begin{eqnarray}\label{eqn:proof-laplace-comparison}\nonumber\langle\langle\Delta_t u, \, u\rangle\rangle  & = & ||d\, u||^2 + ||d_t^{\star}\, u||^2  \\
\nonumber  & = & ||\partial_t u||^2 + ||\bar\partial_t u||^2 + ||\partial_t^{\star}u||^2 + ||\bar\partial_t^{\star}u||^2 \\
\nonumber  & = & \langle\langle\Delta_t' u, \, u\rangle\rangle + \langle\langle\Delta_t'' u, \, u\rangle\rangle \\
           & \geq & \langle\langle\Delta_t' u, \, u\rangle\rangle\geq 0,\end{eqnarray}

\noindent since $\langle\langle\Delta_t' u, \, u\rangle\rangle = ||\partial_t u||^2 + ||\partial_t^{\star}u||^2\geq 0$ and $\langle\langle\Delta_t'' u, \, u\rangle\rangle = ||\bar\partial_t u||^2 + ||\bar\partial_t^{\star}u||^2\geq 0$, while the assumption that $u$ is not $\Delta_t''$-harmonic amounts to $\langle\langle\Delta_t'' u, \, u\rangle\rangle >0$. The equality between the top two lines follows from $d\, u=\partial_t u + \bar\partial_t u$ and the pure-type forms $\partial_t u$ and $\bar\partial_t u$ of distinct types $(r+1, \, s)$, resp. $(r, \, s+1)$, being orthogonal. Thus $||d\, u||^2=||\partial_t u||^2 + ||\bar\partial_t u||^2$ and the adjoints satisfy the analogous identity $||d_t^{\star}\, u||^2=||\partial_t^{\star} u||^2 + ||\bar\partial_t^{\star} u||^2$ for the same reasons.

 Thus it follows from (\ref{eqn:half-lap-norms}) and (\ref{eqn:proof-laplace-comparison}) that 

\begin{equation}\label{eqn:xi-Omega-ref}||\xi_{t, \, (p+1)}||\leq ||\Omega_{t, \, (p)}^{n-1, \, n-1}||.\end{equation} 

 Now $||\xi_{t, \, (p+1)}^{n-1, \, n-1}||\leq ||\xi_{t, \, (p+1)}||$ by mutual orthogonality of the pure-type components of $\xi_{t, \, (p+1)}$. Similarly $||\xi_{t, \, (p)}^{n-1, \, n-1}||\leq ||\xi_{t, \, (p)}||$, while 

\begin{equation}\label{eqn:L2norm-min-del}||\Omega_{t, \, (p)}^{n-1, \, n-1}||\leq ||\xi_{t, \, (p)}^{n-1, \, n-1}||\end{equation} 

\noindent by $L^2$-norm minimality of $\Omega_{t, \, (p)}^{n-1, \, n-1}$ among the solutions of the equation $\partial_t\Omega_{t, \, (p)}^{n-1, \, n-1}=\partial_t\xi_{t, \, (p)}^{n-1, \, n-1}$ (cf. (\ref{eqn:iterative-pot-stepp})). Thus we get

$$||\xi_{t, \, (p+1)}^{n-1, \, n-1}||\leq ||\xi_{t, \, (p+1)}|| \leq ||\Omega_{t, \, (p)}^{n-1, \, n-1}||\leq ||\xi_{t, \, (p)}^{n-1, \, n-1}||\leq ||\xi_{t, \, (p)}||.$$

\noindent This sequence of inequalities contains (\ref{eqn:decrease}).  \hfill $\Box$

\vspace{2ex}

 Taking our cue from the strict inequality (\ref{eqn:laplace-comparison-strict}), we now notice that inequality (\ref{eqn:decrease}) can be improved in a way that is uniform w.r.t. $t\in\Delta$ if the relevant forms $\partial_t\xi_{t, \, (p)}^{n-1, \, n-1}$ avoid the harmonic spaces $\ker\Delta_t''$ for all $t\in\Delta$
(including $t=0$). This is not possible, however, if the non-$\Delta_t''$-harmonicity assumption is only made at $t\neq 0$.

\begin{Obs}\label{Obs:t-unif} \,\, $(i)$ Let $(u_t)_{t\in\Delta}$ be a family of $J_t$-$(r, \, s)$-forms varying continuously with $t$ (up to $t=0$) such that $u_t\notin\ker\Delta_t''$ for all $t\in\Delta$ (including $t=0$). 

 Then there exists a constant $\varepsilon >0$ independent of $t\in\Delta$ such that

\begin{equation}\label{eqn:epsilont-laplace-comparison}\langle\langle\Delta_t u_t, \, u_t\rangle\rangle \geq (1+\varepsilon)\, \langle\langle\Delta_t' u_t, \, u_t\rangle\rangle \hspace{2ex}  \mbox{for all} \hspace{2ex} t\in\Delta,\end{equation}

\noindent after possibly shrinking the base $\Delta$ about $0$. \\

\noindent $(ii)$\, In particular, suppose that for a given $p\in\N$ we have

\begin{equation}\label{eqn:unif-hyp}\partial_t\xi_{t, \, (p)}^{n-1, \, n-1}\notin \ker\Delta_t'',  \hspace{3ex}  \mbox{for all} \hspace{2ex} t\in\Delta  \hspace{2ex}  (\mbox{including} \hspace{2ex} t=0).\end{equation}

 Then there exists a constant $\varepsilon >0$ independent of $t\in\Delta$ such that

\begin{equation}\label{eqn:oneplusepsilon-t-est}\langle\langle\Delta_t(\partial_t\xi_{t, \, (p)}^{n-1, \, n-1}), \, \partial_t\xi_{t, \, (p)}^{n-1, \, n-1}\rangle\rangle\geq (1+\varepsilon)\,\langle\langle\Delta_t'(\partial_t\xi_{t, \, (p)}^{n-1, \, n-1}), \, \partial_t\xi_{t, \, (p)}^{n-1, \, n-1}\rangle\rangle, \hspace{2ex} t\in\Delta,\end{equation}

\noindent after possibly shrinking the base $\Delta$ about $0$.

 Implicitly, if hypothesis (\ref{eqn:unif-hyp}) is satisfied for a given $p\in\N$, we get

\begin{equation}\label{eqn:epsilon-est-xi-p}||\xi_{t, \, (p+1)}^{n-1, \, n-1}||\leq \frac{1}{\sqrt{1+\varepsilon}}\, ||\xi_{t, \, (p)}^{n-1, \, n-1}||, \hspace{3ex} t\in\Delta,\end{equation}

\noindent for a certain $\varepsilon>0$ independent of $t\in\Delta$.

\end{Obs}

\noindent {\it Proof.} Since $u_0\notin\ker\Delta_0''$, inequality (\ref{eqn:laplace-comparison-strict}) applies to $\Delta_0$, $\Delta_0'$ and $u_0$ to give $\langle\langle\Delta_0 u_0, \, u_0\rangle\rangle > \langle\langle\Delta_0' u_0, \, u_0\rangle\rangle.$ Thus there exists a constant $\varepsilon >0$ such that this inequality strengthens to

\begin{equation}\label{eqn:epsilon-laplace-comparison}\langle\langle\Delta_0 u_0, \, u_0\rangle\rangle > (1+\varepsilon)\, \langle\langle\Delta_0' u_0, \, u_0\rangle\rangle.\end{equation}
 
 Now $(\Delta_t)_{t\in\Delta}$ and $(\Delta_t')_{t\in\Delta}$ are $C^{\infty}$ families of operators since they are defined by metrics $(\gamma_t)_{t\in\Delta}$ that vary in a $C^{\infty}$ way with $t\in\Delta$ (up to $t=0$). As a result, $\langle\langle\Delta_t u_t, \, u_t\rangle\rangle$ and $\langle\langle\Delta_t' u_t, \, u_t\rangle\rangle$ both vary continuously with $t\in\Delta$. By continuity, shrinking $\Delta$ about $0$ if necessary, (\ref{eqn:epsilon-laplace-comparison}) extends to a small neighbourhood of $0$ in $\Delta$ to give (\ref{eqn:epsilont-laplace-comparison}) and prove part $(i)$.

 The first statement of part $(ii)$ is an immediate consequence of part $(i)$. As for the second statement of part $(ii)$, it follows from (\ref{eqn:oneplusepsilon-t-est}), (\ref{eqn:re-min-form-p}), (\ref{eqn:delta'-halfnorm}) and (\ref{eqn:delta-halfnorm}) that

$$||\xi_{t, \, (p+1)}||\leq\frac{1}{\sqrt{1+\varepsilon}}\, ||\Omega_{t, \, (p)}^{n-1, \, n-1}||,$$

\noindent while the easy comparison arguments given at the end of the proof of Lemma \ref{Lem:decrease} further give the uniform estimate (\ref{eqn:epsilon-est-xi-p}). The proof is complete.   \hfill $\Box$

\vspace{3ex}

 The forms $\xi_{t, \, (p)}^{n-1, \, n-1}$ produced iteratively in Lemma \ref{Lem:p-iterations} may appear at first glance as the right substitute for the previous forms $\xi_t^{n-1, \, n-1}$ if $p\gg 1$. However, the $L^2$-norm of $\xi_{t, \, (p)}^{n-1, \, n-1}$ need not be small uniformly w.r.t. $t\in\Delta$ and the number $p\gg 1$ of iterations due to the uncontrollable behaviour of $\partial_t\xi_{t, \, (p)}^{n-1, \, n-1}$ from which $\xi_{t, \, (p+1)}^{n-1, \, n-1}$ is constructed by solving equations (\ref{eqn:iterative-pot-stepp}). Indeed, $\partial_t\xi_{t, \, (p)}^{n-1, \, n-1}$ cannot be guaranteed to satisfy hypothesis (\ref{eqn:unif-hyp}) for all $t\in\Delta$ and all $p\in\N$. Consequently, estimate (\ref{eqn:epsilon-est-xi-p}) need not hold at all, let alone with a constant $\varepsilon>0$ independent of both $t\in\Delta$ and $p\in\N$. Even in the favourable case where $\partial_t\xi_{t, \, (p)}^{n-1, \, n-1}\notin\ker\Delta_t''$ for all $t$ and $p$, the $\varepsilon$ of (\ref{eqn:oneplusepsilon-t-est}) cannot not be guaranteed to be independent of $p\in\N$ since $\partial_t\xi_{t, \, (p)}^{n-1, \, n-1}$ may come arbitrarily close to $\ker\Delta_t''$ as $p\rightarrow +\infty$.

 In other words, we cannot guarantee that inequality (\ref{eqn:decrease}) does not become an identity for $p\gg 1$ or that the decrease of the $L^2$-norms $||\xi_{t, \, (p)}^{n-1, \, n-1}||$, should it occur, is uniform w.r.t. $t$ and $p$ as $p\rightarrow +\infty$. A further modification is needed to achieve uniformity in the $L^2$-estimates: the forms $\xi_{t, \, (p)}$ will be replaced by new inductively constructed forms $\widetilde\xi_{t, \, (p)}$ obtained in the following way. If $\widetilde\xi_{t, \, (p)}$ has been constructed at step $p$ of the inductive procedure that will be described below, $\partial_t\widetilde\xi_{t, \, (p)}^{n-1, \, n-1}$ will be altered to $\partial_t(\widetilde\xi_{t, \, (p)}^{n-1, \, n-1}+\nu_{t, \, (p)}^{n-1, \, n-1})$ (for a suitably chosen form $\nu_{t, \, (p)}^{n-1, \, n-1}$ that will force $\partial_t(\widetilde\xi_{t, \, (p)}^{n-1, \, n-1}+\nu_{t, \, (p)}^{n-1, \, n-1})$ to satisfy analogues of hypothesis (\ref{eqn:unif-hyp}) and of estimate (\ref{eqn:oneplusepsilon-t-est})) before solving equations analogous to (\ref{eqn:iterative-pot-stepp}) and running step $(p+1)$ of the inductive procedure that will produce the next form $\widetilde\xi_{t, \, (p+1)}$. Thus we will ``push'' $\partial_t\widetilde\xi_{t, \, (p)}^{n-1, \, n-1}$ away from $\ker\Delta_t''$ by adding some auxiliary form $\partial_t\nu_{t, \, (p)}^{n-1, \, n-1}$ changing with $p$. We stress that the auxiliary form must be changed at every step $p$ to shift $\partial_t\widetilde\xi_{t, \, (p)}^{n-1, \, n-1}$ beyond a uniform distance from $\ker\Delta_t''$. There is no ``universal'' choice of auxiliary form that would suit every $p$. The details are spelt out in the next sections.

 \subsection{The new inductive construction}\label{subsection:step1_inductive-construction}

 \noindent {\bf Step $1$.} By (\ref{eqn:half-conclusion3}) of Lemma \ref{Lem:summingup} we get

\begin{equation}\label{eqn:half-conclusion3bis}\partial_t\gamma_t^{n-1}=\partial_t\xi_t^{n-1, \, n-1}+\bar\partial_t(\xi_t^{n, \, n-2}+w_t), \hspace{3ex} t\in\Delta.\end{equation}

 Let $(\eta_t)_{t\in\Delta}$ be a smooth family of $J_t$-$(n, \, n-1)$-forms (the {\it auxiliary forms at step $1$}) satisfying the following three conditions ($\star$): \\

\noindent $(a)$\,\, $\eta_t=\partial_t\nu_t^{n-1, \, n-1}=\bar\partial_t\vartheta_t^{n, \, n-2}$ for all $t\in\Delta$ and for continuous families of forms $(\nu_t^{n-1, \, n-1})_{t\in\Delta}$, $(\vartheta_t^{n, \, n-2})_{t\in\Delta}$ of the shown types; \\

\noindent $(b)$\,\, $||\xi_t^{n-1, \, n-1} + \nu_t^{n-1, \, n-1}||\leq ||\xi_t||, \hspace{3ex} t\in\Delta;$ \\

\noindent $(c)$\,\, for all $t\in\Delta$ and for some $\varepsilon_0>0$ independent of $t$ we have \\

$\displaystyle\frac{\langle\langle\Delta_t''(\partial_t\xi_t^{n-1, \, n-1} + \partial_t\nu_t^{n-1, \, n-1}), \, \partial_t\xi_t^{n-1, \, n-1} + \partial_t\nu_t^{n-1, \, n-1}\rangle\rangle}{\langle\langle\Delta_t'(\partial_t\xi_t^{n-1, \, n-1} + \partial_t\nu_t^{n-1, \, n-1}), \, \partial_t\xi_t^{n-1, \, n-1} + \partial_t\nu_t^{n-1, \, n-1}\rangle\rangle} \geq \varepsilon_0>0,$ \\

\noindent with the convention that if the denominator vanishes, any $\varepsilon_0>0$ will do.

\vspace{2ex}

\noindent Now using $(a)$, (\ref{eqn:half-conclusion3bis}) becomes:

\begin{equation}\label{eqn:half-conclusion3bisbis}\partial_t\gamma_t^{n-1}=\partial_t(\xi_t^{n-1, \, n-1}+\nu_t^{n-1, \, n-1}) + \bar\partial_t(\xi_t^{n, \, n-2}+w_t-\vartheta_t^{n, \, n-2}), \hspace{3ex} t\in\Delta.\end{equation}

 Let $\widetilde{\Omega}_t^{n-1, \, n-1}$ and $\widetilde{\xi}_{t, \, (1)}$ be the $\partial_t$-potential and respectively the $d$-potential of minimal $L^2$-norms of $\partial_t(\xi_t^{n-1, \, n-1} + \nu_t^{n-1, \, n-1}):$

\begin{equation}\label{eqn:re-iterative-pot-1}\partial_t(\xi_t^{n-1, \, n-1} + \nu_t^{n-1, \, n-1}) = \partial_t\widetilde{\Omega}_t^{n-1, \, n-1} = d\, \widetilde{\xi}_{t, \, (1)}, \hspace{2ex} t\in\Delta.\end{equation} 

\noindent Notice that, since $d\, \widetilde{\xi}_{t, \, (1)}$ is of pure type $(n, \, n-1)$, we must have 

$$d\, \widetilde{\xi}_{t, \, (1)}=\partial_t\widetilde{\xi}_{t, \, (1)}^{n-1, \, n-1} + \bar\partial_t\widetilde{\xi}_{t, \, (1)}^{n, \, n-2}, \hspace{2ex} t\in\Delta.$$

 Using this and (\ref{eqn:re-iterative-pot-1}), (\ref{eqn:half-conclusion3bisbis}) reads: 

\begin{equation}\label{eqn:half-conclusion3bisbisfinal}\partial_t\gamma_t^{n-1}=\partial_t\widetilde{\xi}_{t, \, (1)}^{n-1, \, n-1} + \bar\partial_t(\widetilde{\xi}_{t, \, (1)}^{n, \, n-2} + \xi_t^{n, \, n-2}+w_t-\vartheta_t^{n, \, n-2}), \hspace{3ex} t\in\Delta.\end{equation}

\noindent {\bf Step $p+1$.} Suppose that Step $p$ has been performed and has produced the following decomposition for all $t\in\Delta$:

\begin{eqnarray}\label{eqn:half-conclusion3bisbisfinalp}\nonumber\partial_t\gamma_t^{n-1}=\partial_t\widetilde{\xi}_{t, \, (p)}^{n-1, \, n-1} & + & \bar\partial_t(\widetilde{\xi}_{t, \, (p)}^{n, \, n-2} + \dots + \widetilde{\xi}_{t, \, (1)}^{n, \, n-2} + \xi_t^{n, \, n-2}+w_t  \\
  & - & \vartheta_t^{n, \, n-2} - \vartheta_{t, \, (1)}^{n, \, n-2} - \dots - \vartheta_{t, \, (p-1)}^{n, \, n-2}).\end{eqnarray}

 Let $(\eta_{t, \, (p)})_{t\in\Delta}$ be a smooth family of $J_t$-$(n, \, n-1)$-forms (the {\it auxiliary forms at step $p+1$}) satisfying the following three conditions ($\star_p$): \\

\noindent $(a)$\,\, $\eta_{t, \, (p)}=\partial_t\nu_{t, \, (p)}^{n-1, \, n-1}=\bar\partial_t\vartheta_{t, \, (p)}^{n, \, n-2}$ for all $t\in\Delta$ and for continuous families of forms $(\nu_{t, \, (p)}^{n-1, \, n-1})_{t\in\Delta}$, $(\vartheta_{t, \, (p)}^{n, \, n-2})_{t\in\Delta}$ of the shown types; \\

\noindent $(b)$\,\, $||\widetilde{\xi}_{t, \, (p)}^{n-1, \, n-1} + \nu_{t, (p)}^{n-1, \, n-1}||\leq ||\widetilde{\xi}_{t, \, (p)}||, \hspace{3ex} t\in\Delta;$ \\

\noindent $(c)$\,\, for all $t\in\Delta$ and for some $\varepsilon_0>0$ independent of $t$ and of $p\in\N$ we have \\

$\displaystyle\frac{\langle\langle\Delta_t''(\partial_t\widetilde{\xi}_{t, \, (p)}^{n-1, \, n-1} + \partial_t\nu_{t, \, (p)}^{n-1, \, n-1}), \, \partial_t\widetilde{\xi}_{t, \, (p)}^{n-1, \, n-1} + \partial_t\nu_{t, \, (p)}^{n-1, \, n-1}\rangle\rangle}{\langle\langle\Delta_t'(\partial_t\widetilde{\xi}_{t, \, (p)}^{n-1, \, n-1} + \partial_t\nu_{t, \, (p)}^{n-1, \, n-1}), \, \partial_t\widetilde{\xi}_{t, \, (p)}^{n-1, \, n-1} + \partial_t\nu_{t, \, (p)}^{n-1, \, n-1}\rangle\rangle} \geq \varepsilon_0>0,$ \\

\noindent with the convention that if the denominator vanishes, any $\varepsilon_0>0$ will do.

\vspace{2ex}

\noindent Now using $(a)$, (\ref{eqn:half-conclusion3bisbisfinalp}) becomes for all $t\in\Delta$:

\begin{eqnarray}\label{eqn:half-conclusion3bisbisp+1}\nonumber\partial_t\gamma_t^{n-1}=\partial_t(\widetilde{\xi}_{t, \, (p)}^{n-1, \, n-1}+\nu_{t, \, (p)}^{n-1, \, n-1}) & + & \bar\partial_t(\widetilde\xi_{t, \, (p)}^{n, \, n-2}+ \dots + \widetilde\xi_{t, \, (1)}^{n, \, n-2}+ \xi_t^{n, \, n-2} + w_t \\
   & - & \vartheta_t^{n, \, n-2} - \vartheta_{t, \, (1)}^{n, \, n-2} - \dots - \vartheta_{t, \, (p)}^{n, \, n-2}).\end{eqnarray}

 Let $\widetilde{\Omega}_{t, \, (p)}^{n-1, \, n-1}$ and $\widetilde{\xi}_{t, \, (p+1)}$ be the $\partial_t$-potential and respectively the $d$-potential of minimal $L^2$-norms of $\partial_t(\widetilde{\xi}_{t, \, (p)}^{n-1, \, n-1} + \nu_{t, \, (p)}^{n-1, \, n-1}):$

\begin{equation}\label{eqn:re-iterative-pot}\partial_t(\widetilde{\xi}_{t, \, (p)}^{n-1, \, n-1} + \nu_{t, \, (p)}^{n-1, \, n-1}) = \partial_t\widetilde{\Omega}_{t, \, (p)}^{n-1, \, n-1} = d\, \widetilde{\xi}_{t, \, (p+1)}, \hspace{2ex} t\in\Delta.\end{equation} 

\noindent Notice that, since $d\, \widetilde{\xi}_{t, \, (p+1)}$ is of pure type $(n, \, n-1)$, we must have 

$$d\, \widetilde{\xi}_{t, \, (p+1)}=\partial_t\widetilde{\xi}_{t, \, (p+1)}^{n-1, \, n-1} + \bar\partial_t\widetilde{\xi}_{t, \, (p+1)}^{n, \, n-2}, \hspace{2ex} t\in\Delta.$$

\noindent Using this and (\ref{eqn:re-iterative-pot}), (\ref{eqn:half-conclusion3bisbisp+1}) reads for all $t\in\Delta$:

\begin{eqnarray}\label{eqn:half-conclusion3bisbisfinal}\nonumber\partial_t\gamma_t^{n-1}=\partial_t\widetilde{\xi}_{t, \, (p+1)}^{n-1, \, n-1} & + & \bar\partial_t(\widetilde{\xi}_{t, \, (p+1)}^{n, \, n-2} + \dots + \widetilde{\xi}_{t, \, (1)}^{n, \, n-2} + \xi_t^{n, \, n-2} + w_t \\
 & - & \vartheta_t^{n, \, n-2} - \vartheta_{t, \, (1)}^{n, \, n-2} - \dots - \vartheta_{t, \, (p)}^{n, \, n-2}),\end{eqnarray}

\noindent completing the inductive construction of the families $(\widetilde{\xi}_{t, \, (p)}^{n-1, \, n-1})_{t\in\Delta}$, $p\in\N$. 

\vspace{2ex}

 {\bf Summing up:} if we set $\widetilde{\xi}_{t, \, (0)}: = \xi_t$ and $\widetilde{\xi}_{t, \, (0)}^{n-1, \, n-1}: = \xi_t^{n-1, \, n-1}$ as well as $\widetilde{\Omega}_{t, \, (0)}^{n-1, \, n-1}:= \widetilde{\Omega}_t^{n-1, \, n-1}$ and $\nu_{t, \, (0)}^{n-1, \, n-1}:= \nu_t^{n-1, \, n-1}$, we get continuous families of forms $(\widetilde{\xi}_{t, \, (p)}^{n-1, \, n-1})_{t\in\Delta}$ and $(\widetilde{\Omega}_{t, \, (p)}^{n-1, \, n-1})_{t\in\Delta}$ for each $p\in\N$.

\begin{Com}\label{Com:zero-aux-forms} {\rm It is clear that the forms $\eta_{t, \, (p)}=0$ with $\nu_{t, \, (p)}^{n-1, \, n-1}=0$ and $\vartheta_{t, \, (p)}^{n, \, n-2}=0$ for all $t\in\Delta$ trivially satisfy conditions $(a)$ and $(b)$ of $(\star_p)$, while they need not satisfy condition $(c)$. Indeed, we have $||\widetilde{\xi}_{t, \, (p)}^{n-1, \, n-1}|| \leq ||\widetilde{\xi}_{t, \, (p)}||$ (hence $(b)$ for $\nu_{t, \, (p)}^{n-1, \, n-1}=0$) since the former form is the $(n-1, \, n-1)$-component of the latter and forms of distinct pure types are orthogonal. So, in general, the choices of these auxiliary forms are non-trivial. However, the trivial choice of identically zero auxiliary forms will do if it happens to satisfy $(c)$ (see (\ref{eqn:already-c}) below).}

\end{Com}

\section{Proof of the existence of auxiliary forms}\label{section:existence_auxiliary-forms}

 We now spell out the argument accounting for the existence of smooth families of forms $(\eta_{t, \, (p)})_{t\in\Delta}$ satisfying conditions $(\star_p)$ for all $p\in\N$. The spectra of $\Delta_t'$ and $\Delta_t''$ acting on $(n, \, n-1)$-forms satisfy inclusions:

\begin{equation}\label{eqn:spec-not'-''}\mbox{Spec}\,\Delta_t'\subset [0, \, \varepsilon_t']\cup [\varepsilon', \, +\infty), \hspace{2ex} \mbox{Spec}\,\Delta_t''\subset [0, \, \varepsilon_t'']\cup [\varepsilon'', \, +\infty), \hspace{2ex} t\in\Delta,\end{equation}

\noindent where $\varepsilon', \varepsilon''>0$ are independent of $t$, while $\varepsilon_t', \varepsilon_t''\rightarrow 0$ as $t\rightarrow 0$. (Thus $\varepsilon_0'= \varepsilon_0''=0$.) Since the eigenspaces of $\Delta_t'$ and of $\Delta_t''$ are finite-dimensional and since there are at most finitely many eigenvalues of $\Delta_t'$ below $\varepsilon'$ and of $\Delta_t''$ below $\varepsilon''$, each of the vector spaces $\oplus_{\mu\leq\varepsilon'_t}E_{\Delta_t'}^{n, \, n-1}(\mu)$ and $\oplus_{\lambda\leq\varepsilon''_t}E_{\Delta_t''}^{n, \, n-1}(\lambda)$ (which are the obstruction to what we are striving to achieve) has finite dimension. Hence their respective orthogonal complements $\oplus_{\mu\geq\varepsilon'}E_{\Delta_t'}^{n, \, n-1}(\mu)$ and $\oplus_{\lambda\geq\varepsilon''}E_{\Delta_t''}^{n, \, n-1}(\lambda)$ in the infinite-dimensional vector space $C^{\infty}_{n, \, n-1}(X_t, \, \C)$ have both infinite dimension and so has their intersection with the infinite-dimensional subspace $\mbox{Im}\,\bar\partial_t$, i.e.

$${\cal E}_t^{n, \, n-1}:=\bigoplus\limits_{\mu\geq\varepsilon'}E_{\Delta_t'}^{n, \, n-1}(\mu)\cap\bigoplus\limits_{\lambda\geq\varepsilon''}E_{\Delta_t''}^{n, \, n-1}(\lambda)\cap\mbox{Im}\,\bar\partial_t\subset C^{\infty}_{n, \, n-1}(X_t, \, \C), \hspace{2ex} t\in\Delta,$$

\noindent has infinite dimension. The infinite dimensionality of ${\cal E}_t^{n, \, n-1}$ will play a crucial role in the sequel\!: the auxiliary forms $\eta_{t, \, (p)}$ will be chosen in ${\cal E}_t^{n, \, n-1}$ and having plenty of ``room for choice'' will be a key factor. Moreover, $\Delta\ni t\mapsto{\cal E}_t^{n, \, n-1}$ defines an infinite-rank $C^{\infty}$-subbundle of $\Delta\ni t\mapsto C^{\infty}_{n, \, n-1}(X_t, \, \C)$. Notice the inclusion

\begin{equation}\label{eqn:del-incl}{\cal E}_t^{n, \, n-1}\subset\mbox{Im}\,\partial_t, \hspace{3ex} t\in\Delta.\end{equation}

\noindent Indeed, being of type $(n, \, n-1)$, every form $\eta_t\in{\cal E}_t^{n, \, n-1}$ is trivially $\partial_t$-closed, hence also $d$-closed since the $\bar\partial_t$-exactness assumption is implicit in the definition of ${\cal E}_t^{n, \, n-1}$. Then $\eta_t$ is $\partial_t$-exact for all $t\neq 0$ by the $\partial\bar\partial$-lemma. Since any $\eta_t\in{\cal E}_t^{n, \, n-1}$ avoids the {\it small} eigenvalues of $\Delta_t'$ by definition of ${\cal E}_t^{n, \, n-1}$, it follows that $\eta_0$ must be again $\partial_0$-exact if $\eta_0$ stands in a $C^{\infty}$ family $(\eta_t)_{t\in\Delta}$ with $\eta_t\in{\cal E}_t^{n, \, n-1}$ for all $t\in\Delta$.

 Now fix $p\in\N$ and suppose that the induction has been performed up to {\it Step} $p$. In particular, the forms $(\widetilde{\xi}_{t, \, (p)})_{t\in\Delta}$ have already been constructed. To run {\it Step} $(p+1)$, we have to show the existence of auxiliary forms $(\eta_{t, \, (p)})_{t\in\Delta}$ adapted to the pre-existing forms $(\widetilde{\xi}_{t, \, (p)}^{n-1, \, n-1})_{t\in\Delta}$ by satisfying conditions $(\star_p)$. To start with, pick any smooth family $(\eta_{t, \, (p)})_{t\in\Delta}$ of non-zero $J_t$-$(n, \, n-1)$-forms such that 

\begin{equation}\label{eqn:choiceextra}\eta_{t, \, (p)}=\partial_t\nu_{t, \, (p)}^{n-1, \, n-1}=\bar\partial_t\vartheta_{t, \, (p)}^{n, \, n-2}\in{\cal E}_t^{n, \, n-1}, \hspace{2ex} t\in\Delta,\end{equation}

\noindent where the families of {\it minimal $L^2$-norm} $\partial_t$-potentials $(\nu_{t, \, (p)}^{n-1, \, n-1})_{t\in\Delta}$ and $\bar\partial_t$-potentials $(\vartheta_{t, \, (p)}^{n, \, n-2})_{t\in\Delta}$ vary continuously with $t\in\Delta$ (up to $t=0$).  We thus satisfy requirement $(a)$ in the infinite-rank vector bundle $\Delta\ni t\mapsto {\cal E}_t^{n, \, n-1}\subset C^{\infty}_{n, \, n-1}(X_t, \, \C)$. We have yet to satisfy the requirements $(b)$ and $(c)$. 

\vspace{2ex}

 Fo every $t\in\Delta$, consider the map

\begin{equation}\label{eqn:defS_t}{\cal E}_t^{n, \, n-1}\ni\eta_{t, \, (p)}\stackrel{S_t}{\longmapsto} \nu_{t, \, (p)}^{n-1, \, n-1}\in\mbox{Im}\,\partial_t^{\star}\subset C_{n-1, \, n-1}^{\infty}(X_t, \, \C)\end{equation}

\noindent which associates with every $\eta_{t, \, (p)}\in {\cal E}_t^{n, \, n-1}$ its $\partial_t$-potential $\nu_{t, \, (p)}^{n-1, \, n-1}$ of {\it minimal $L^2$-norm} (i.e. the unique $\partial_t$-potential that lies in $\mbox{Im}\,\partial_t^{\star}$). Since $\mbox{Im}\, S_t\subset\mbox{Im}\,\partial_t^{\star}$, we have 

\begin{equation}\label{eqn:S-t-orthog}\mbox{Im}\, S_t\perp\ker\partial_t.\end{equation}

\noindent It is clear that the map $S_t$ is linear (because $\nu_{t, \, (p)}^{n-1, \, n-1}=\Delta_t^{'-1}\partial_t^{\star}\eta_{t, \, (p)}$ while $\Delta_t^{'-1}$ and $\partial_t^{\star}$ are linear operators) and injective (because $\partial_t\nu_{t, \, (p)}^{n-1, \, n-1}=\eta_{t, \, (p)}$). Hence $\mbox{Im}\,S_t$ is an infinite-dimensional vector subspace of $\mbox{Im}\,\partial_t^{\star}$.

  Meanwhile, for every $t\in\Delta$ and every $p\in\N$, let

 $${\cal U}_{t, \, (p)}:=\overline{B}\bigg(-\widetilde{\xi}_{t, \, (p)}^{n-1, \, n-1}, \, ||\widetilde{\xi}_{t, \, (p)}||\bigg)\subset C^{\infty}_{n-1, \, n-1}(X_t, \, \C)$$

\noindent be the closed ball (w.r.t. $L^2$-norm) centred at $-\widetilde{\xi}_{t, \, (p)}^{n-1, \, n-1}$ and of radius $||\widetilde{\xi}_{t, \, (p)}||$ in $C^{\infty}_{n-1, \, n-1}(X_t, \, \C)$. Clearly, $0\in{\cal U}_{t, \, (p)}$. Condition $(b)$ of $(\star_p)$ translates to

\begin{equation}\label{eqn:nu-ball-belong}\nu_{t, \, (p)}^{n-1, \, n-1}\in{\cal U}_{t, \, (p)}, \hspace{3ex}  t\in\Delta,\end{equation}

\noindent so any form

\begin{equation}\label{eqn:ball-inter-ImS}\nu_{t, \, (p)}^{n-1, \, n-1}\in\mbox{Im}\,S_t\cap{\cal U}_{t, \, (p)}, \hspace{3ex}  t\in\Delta,\end{equation}

\noindent automatically satisfies conditions $(a)$ (after setting $\eta_{t, \, (p)}\!\!:=\partial_t\nu_{t, \, (p)}^{n-1, \, n-1}$) and $(b)$ of $(\star_p)$. Note that unless the form $\widetilde{\xi}_{t, \, (p)}^{n-1, \, n-1}$ (given by the induction hypothesis) already satisfies the condition

\begin{equation}\label{eqn:already-c}\frac{\langle\langle\Delta_t''(\partial_t\widetilde{\xi}_{t, \, (p)}^{n-1, \, n-1}), \, \partial_t\widetilde{\xi}_{t, \, (p)}^{n-1, \, n-1}\rangle\rangle}{\langle\langle\Delta_t'(\partial_t\widetilde{\xi}_{t, \, (p)}^{n-1, \, n-1}), \, \partial_t\widetilde{\xi}_{t, \, (p)}^{n-1, \, n-1}\rangle\rangle} \geq \varepsilon_0>0,\end{equation}

\noindent for the uniform $\varepsilon_0$ obtained from the previous induction steps $1, \dots , p$, the auxiliary form $\nu_{t, \, (p)}^{n-1, \, n-1}$ that we are now trying to construct cannot be chosen to be the zero form. Thus, unless $\widetilde{\xi}_{t, \, (p)}^{n-1, \, n-1}$ satisfies (\ref{eqn:already-c}), we must show that

\begin{equation}\label{eqn:ball-inter-ImS-notzero}\mbox{Im}\,S_t\cap{\cal U}_{t, \, (p)}\supsetneq\{0\}, \hspace{3ex}  t\in\Delta.\end{equation}

\noindent If we can manage to achieve (\ref{eqn:ball-inter-ImS-notzero}), we will choose $0\neq\nu_{t, \, (p)}^{n-1, \, n-1}\in\mbox{Im}\,S_t\cap{\cal U}_{t, \, (p)}$ (cf. (\ref{eqn:ball-inter-ImS})) in a family varying in a continuous way with $t\in\Delta$ and will set $\eta_{t, \, (p)}:= \partial_t\nu_{t, \, (p)}^{n-1, \, n-1}$ for every $t\in\Delta$. Property (\ref{eqn:choiceextra}) will then be satisfied and so will be $(a)$ and $(b)$ of $(\star_p)$.

 The discussion of the possibility of enforcing the choice (\ref{eqn:ball-inter-ImS}) falls into two cases that we now analyse. \\

\noindent {\it Case $1$: if $||\widetilde{\xi}_{t, \, (p)}^{n-1, \, n-1}|| < ||\widetilde{\xi}_{t, \, (p)}||$}, then the origin $0$ of $C^{\infty}_{n-1, \, n-1}(X_t, \, \C)$ lies in the interior of the ball ${\cal U}_{t, \, (p)}$, so the vector subspace $\mbox{Im}\, S_t$ meets the interior of ${\cal U}_{t, \, (p)}$. Hence (\ref{eqn:ball-inter-ImS-notzero}) is guaranteed and we can choose $\nu_{t, \, (p)}^{n-1, \, n-1}\neq 0$ to satisfy  (\ref{eqn:ball-inter-ImS}). Conditions $(a)$ and $(b)$ of $(\star_p)$ are thus simultaneously fulfilled as explained above.

\vspace{2ex}

\noindent {\it Case $2$: if $||\widetilde{\xi}_{t, \, (p)}^{n-1, \, n-1}|| = ||\widetilde{\xi}_{t, \, (p)}||$}, then $\widetilde{\xi}_{t, \, (p)}=\widetilde{\xi}_{t, \, (p)}^{n-1, \, n-1}$, hence $\widetilde{\xi}_{t, \, (p)}$ is of pure type $(n-1, \, n-1)$. (Recall that, in general, $||\widetilde{\xi}_{t, \, (p)}||^2 = ||\widetilde{\xi}_{t, \, (p)}^{n, \, n-2}||^2 + ||\widetilde{\xi}_{t, \, (p)}^{n-1, \, n-1}||^2 + ||\widetilde{\xi}_{t, \, (p)}^{n-2, \, n}||^2$ by mutual orthogonality of the pure-type components of a given form.) In this case the zero form $0$ lies on the boundary of the ball ${\cal U}_{t, \, (p)}$. 

 Let  $H_{t, \, (p)}$ denote the hyperplane of $C^{\infty}_{n-1, \, n-1}(X_t, \, \C)$ that is orthogonal to the vector $\widetilde{\xi}_{t, \, (p)}^{n-1, \, n-1}$ at $0$. If the inclusion

\begin{equation}\label{eqn:ImS-incl-hyper}\mbox{Im}\, S_t\subset H_{t, \, (p)}\end{equation} 

\noindent does not hold, then $\mbox{Im}\, S_t$ meets the interior of the ball ${\cal U}_{t, \, (p)}$, (\ref{eqn:ball-inter-ImS-notzero}) holds, we can choose $\nu_{t, \, (p)}^{n-1, \, n-1}\neq 0$ to satisfy (\ref{eqn:ball-inter-ImS}) and we can proceed as in {\it Case $1$}. 

 However, if the inclusion (\ref{eqn:ImS-incl-hyper}) happens to hold, then $\mbox{Im}\, S_t$ does not meet the interior of ${\cal U}_{t, \, (p)}$ and $\mbox{Im}\,S_t\cap{\cal U}_{t, \, (p)}=\{0\}$. Thus (\ref{eqn:ball-inter-ImS-notzero}) does not hold. Meanwhile recall that $\widetilde{\xi}_{t, \, (p)}$ satisfies (by construction) the following induction hypothesis (cf. (\ref{eqn:re-iterative-pot}) with $p-1$ in place of $p$):

\begin{equation}\label{eqn:recall-ind-hyp}\partial_t(\widetilde{\xi}_{t, \, (p-1)}^{n-1, \, n-1} + \nu_{t, \, (p-1)}^{n-1, \, n-1}) = \partial_t\widetilde{\Omega}_{t, \, (p-1)}^{n-1, \, n-1} = d\, \widetilde{\xi}_{t, \, (p)}, \hspace{2ex} t\in\Delta.\end{equation}

\noindent Since $d\, \widetilde{\xi}_{t, \, (p)}= \partial_t\, \widetilde{\xi}_{t, \, (p)} + \bar\partial_t\, \widetilde{\xi}_{t, \, (p)}$ is of pure type $(n, \, n-1)$ and since $\widetilde{\xi}_{t, \, (p)}=\widetilde{\xi}_{t, \, (p)}^{n-1, \, n-1}$ is of pure type $(n-1, \, n-1)$ here, we see that the $(n-1, \, n)$-form $\bar\partial_t\, \widetilde{\xi}_{t, \, (p)}$ must vanish for bidegree reasons. Thus (\ref{eqn:recall-ind-hyp}) yields

\begin{equation}\label{eqn:del-xi-tilde-good}\partial_t(\widetilde{\xi}_{t, \, (p-1)}^{n-1, \, n-1} + \nu_{t, \, (p-1)}^{n-1, \, n-1}) = \partial_t\, \widetilde{\xi}_{t, \, (p)} = \partial_t\,\widetilde{\xi}_{t, \, (p)}^{n-1, \, n-1}, \hspace{2ex} t\in\Delta.\end{equation}

\noindent Now recall that by the induction hypothesis the form $\partial_t(\widetilde{\xi}_{t, \, (p-1)}^{n-1, \, n-1} + \nu_{t, \, (p-1)}^{n-1, \, n-1})$ featuring in the left-hand side of (\ref{eqn:del-xi-tilde-good}) satisfies property $(c)$ of $(\star_{p-1})$ with the uniform $\varepsilon_0>0$ obtained from the previous induction steps $1, \dots , p$. (The auxiliary forms $\nu_{t, \, (p-1)}^{n-1, \, n-1}$ were chosen as such at {\it Step} $p$ of the induction process). Therefore (\ref{eqn:del-xi-tilde-good}) combined with $(c)$ of $(\star_{p-1})$ shows that \\

\noindent $\displaystyle\frac{\langle\langle\Delta_t''(\partial_t\,\widetilde{\xi}_{t, \, (p)}^{n-1, \, n-1}), \, \partial_t\,\widetilde{\xi}_{t, \, (p)}^{n-1, \, n-1}\rangle\rangle}{\langle\langle\Delta_t'(\partial_t\,\widetilde{\xi}_{t, \, (p)}^{n-1, \, n-1}), \, \partial_t\,\widetilde{\xi}_{t, \, (p)}^{n-1, \, n-1}\rangle\rangle} = $

\hspace{7ex} $\displaystyle\frac{\langle\langle\Delta_t''\partial_t(\widetilde{\xi}_{t, \, (p-1)}^{n-1, \, n-1} + \nu_{t, \, (p-1)}^{n-1, \, n-1}), \, \partial_t(\widetilde{\xi}_{t, \, (p-1)}^{n-1, \, n-1} + \nu_{t, \, (p-1)}^{n-1, \, n-1})\rangle\rangle}{\langle\langle\Delta_t'\partial_t(\widetilde{\xi}_{t, \, (p-1)}^{n-1, \, n-1} + \nu_{t, \, (p-1)}^{n-1, \, n-1}), \, \partial_t(\widetilde{\xi}_{t, \, (p-1)}^{n-1, \, n-1} + \nu_{t, \, (p-1)}^{n-1, \, n-1})\rangle\rangle} \geq \varepsilon_0>0,$ \\

\vspace{2ex}

\noindent which means that $\widetilde{\xi}_{t, \, (p)}^{n-1, \, n-1}$ satisfies (\ref{eqn:already-c}). Therefore we can make the trivial choice of auxiliary form
$\nu_{t, \, (p)}^{n-1, \, n-1}$, i.e. we can (and will) choose

$$\nu_{t, \, (p)}^{n-1, \, n-1} = 0\in\mbox{Im}\, S_t\cap{\cal U}_{t, \, (p)} =\{0\}.$$

\noindent This guarantees (\ref{eqn:ball-inter-ImS}), hence $(a)$ and $(b)$ of $(\star_p)$. This also guarantees $(c)$ of $(\star_p)$ thanks to (\ref{eqn:already-c}) (which holds as we have just seen). As explained in Comment \ref{Com:zero-aux-forms}, this choice meets our conditions in this case. (This is the only case where the choice of the {\it zero} form will do.)

\begin{Concl}\label{Concl:sofar} The choice (\ref{eqn:ball-inter-ImS}) can always be enforced and we shall henceforth assume that $\nu_{t, \, (p)}^{n-1, \, n-1}$ has been chosen as in (\ref{eqn:ball-inter-ImS}). This guarantees conditions $(a)$ and $(b)$ of $(\star_p)$. 

 Moreover, in {\it Case} $2$ discussed above, condition $(c)$ is satisfied simultaneously with $(a)$ and $(b)$. It remains to prove that, in {\it Case} $1$ discussed above, $\nu_{t, \, (p)}^{n-1, \, n-1}$ can be chosen as in (\ref{eqn:ball-inter-ImS}) to satisfy furthermore condition $(c)$ of $(\star_p)$.

\end{Concl}

 Let us make the following observation. Since $\nu_{t, \, (p)}^{n-1, \, n-1}$ has been chosen as the minimal $L^2$-norm $\partial_t$-potential of $\eta_{t, \, (p)}$, it satisfies $\nu_{t, \, (p)}^{n-1, \, n-1}\perp\ker\partial_t$ in $C^{\infty}_{n-1, \, n-1}(X_t, \, \C)$ (cf. (\ref{eqn:S-t-orthog})). Thus $\nu_{t, \, (p)}^{n-1, \, n-1}$ cannot have a non-trivial orthogonal projection on any of the eigenspaces $E_{\Delta_t'}^{n-1, \, n-1}(\mu)$ corresponding to eigenvalues $\mu\leq\varepsilon'_t$. Indeed, if $\delta_t\in E_{\Delta_t'}^{n-1, \, n-1}(\mu)\setminus\{0\}$ were such a projection, then $\partial_t\delta_t\in E_{\Delta_t'}^{n, \, n-1}(\mu)\setminus\{0\}$ would play the analogous role for $\eta_{t, \, (p)}=\partial_t\nu_{t, \, (p)}^{n-1, \, n-1}$ in bidegree $(n, \, n-1)$ since $\partial_t$ and $\Delta'_t$ commute. However, the existence of such a component for $\eta_{t, \, (p)}$ is ruled out by (\ref{eqn:choiceextra}) and the definition of ${\cal E}_t^{n, \, n-1}$. Therefore, any form $\nu_{t, \, (p)}^{n-1, \, n-1}\in\mbox{Im}\, S_t$ satisfies 

\begin{equation}\label{eqn:u-ball-in}\nu_{t, \, (p)}^{n-1, \, n-1}\in\bigoplus\limits_{\mu\geq\varepsilon'}E_{\Delta_t'}^{n-1, \, n-1}(\mu),  \hspace{3ex} t\in\Delta.\end{equation}

 We now explain how to choose a form $\nu_{t, \, (p)}^{n-1, \, n-1}$ as in (\ref{eqn:ball-inter-ImS}) that also satisfies requirement $(c)$ of $(\star_p)$ in {\it Case} $1$.

 Condition $(c)$ essentially requires $\partial_t\widetilde\xi_{t, \, (p)}^{n-1, \, n-1} + \partial_t\nu_{t, \, (p)}^{n-1, \, n-1}$ to stay away from $\ker\Delta_t''$ at an $L^2$- distance that is bounded below by a positive constant independent of both $t\in\Delta$ and $p\in\N$ if simultaneously the behaviour of $\partial_t\widetilde\xi_{t, \, (p)}^{n-1, \, n-1} + \partial_t\nu_{t, \, (p)}^{n-1, \, n-1}$ w.r.t. $\Delta_t'$ is kept under control relative to the behaviour w.r.t. $\Delta_t''$. 

  The possibility that $\partial_t\widetilde\xi_{t, \, (p)}^{n-1, \, n-1}$ be $\Delta_t''$-harmonic cannot be ruled out and in this case condition $(c)$ cannot not be fulfilled without correcting $\partial_t\widetilde\xi_{t, \, (p)}^{n-1, \, n-1}$ by non-zero auxiliary forms $\eta_{t, \, (p)}$. Recall that the auxiliary form $\eta_{t, \, (p)}=\partial_t\nu_{t, \, (p)}^{n-1, \, n-1}$ is to be chosen among the forms that satisfy condition (\ref{eqn:choiceextra}). Any such $\eta_{t, \, (p)}$ is $\bar\partial_t$-exact for all $t\in\Delta$ by the choice $(\ref{eqn:choiceextra})$, hence $\eta_{t, \, (p)}$ is orthogonal to $\ker\Delta_t''$ (since $\ker\Delta_t''\perp\mbox{Im}\, \bar\partial_t$). Thus $\eta_{t, \, (p)}=\partial_t\nu_{t, \, (p)}^{n-1, \, n-1}$ is in a good position to ``drive'' $\partial_t\widetilde\xi_{t, \, (p)}^{n-1, \, n-1}$ away from $\ker\Delta_t''$ and ensure that the corrected form $\partial_t\widetilde\xi_{t, \, (p)}^{n-1, \, n-1} + \partial_t\nu_{t, \, (p)}^{n-1, \, n-1}$ satisfies $(c)$.   \\

 The discussion of the choice of a form $\nu_{t, \, (p)}^{n-1, \, n-1}$ as in (\ref{eqn:ball-inter-ImS}) that also satisfies requirement $(c)$ of $(\star_p)$ in {\it Case} $1$ falls into two steps. \\

\noindent $(I)$\,\, {\it Uniformly bounding the numerator of $(c)$ in $(\star_p)$ from below in Case $1$} \\

It is clear that $\langle\langle\Delta_t''(\partial_t\widetilde\xi_{t, \, (p)}^{n-1, \, n-1} + \partial_t\nu_{t, \, (p)}^{n-1, \, n-1}), \, \partial_t\widetilde\xi_{t, \, (p)}^{n-1, \, n-1} + \partial_t\nu_{t, \, (p)}^{n-1, \, n-1}\rangle\rangle$ has a uniform positive lower bound whenever the following three conditions are simultaneously met as $\Delta\ni t\rightarrow 0$ and $p\rightarrow +\infty$: \\

\noindent $(i)$\, the $L^2$-distance from $\partial_t\widetilde\xi_{t, \, (p)}^{n-1, \, n-1} + \partial_t\nu_{t, \, (p)}^{n-1, \, n-1}$ to $\ker\Delta_t''$ does not become arbitrarily small\!; \\

\noindent $(ii)$\, the $L^2$-norm of $\partial_t\widetilde\xi_{t, \, (p)}^{n-1, \, n-1} + \partial_t\nu_{t, \, (p)}^{n-1, \, n-1}$ does not become arbitrarily small\!;\\

\noindent $(iii)$\, $\partial_t\widetilde\xi_{t, \, (p)}^{n-1, \, n-1} + \partial_t\nu_{t, \, (p)}^{n-1, \, n-1}\notin\oplus_{\lambda\leq\varepsilon_t''}E_{\Delta_t''}^{n, \, n-1}(\lambda)$ for $\varepsilon_t''\rightarrow 0$ as $t\rightarrow 0$. \\

 In fact condition $(iii)$ is related to condition $(i)$: if $\partial_t\widetilde\xi_{t, \, (p)}^{n-1, \, n-1} + \partial_t\nu_{t, \, (p)}^{n-1, \, n-1}\in\oplus_{\lambda\leq\varepsilon_t''}E_{\Delta_t''}^{n, \, n-1}(\lambda)$ for $\varepsilon_t''\rightarrow 0$ as $t\rightarrow 0$, then $\partial_0\widetilde\xi_{0, \, (p)}^{n-1, \, n-1} + \partial_0\nu_{0, \, (p)}^{n-1, \, n-1}\in\ker\Delta_0''$ in violation of $(i)$.

\begin{Obs}\label{Obs:nolossofgen} Without loss of generality we may make the following \\

\noindent {\bf Assumption (A1):} \hspace{6ex}   $\partial_0\widetilde\xi_{0, \, (p)}^{n-1, \, n-1}\in\ker\Delta_0''.$ 

\end{Obs}

\noindent {\it Proof.} There are three cases:\\

$(1)$\, if $\partial_t\widetilde\xi_{t, \, (p)}^{n-1, \, n-1}\in\bigoplus\limits_{\lambda\geq\varepsilon''}E^{n, n-1}_{\Delta''_t}(\lambda)$ for all $t\in\Delta$, then $\partial_0\widetilde\xi_{0, \, (p)}^{n-1, \, n-1}\in\mbox{Im}\,\bar\partial_0$. (Indeed, recall that $\partial_t\widetilde\xi_{t, \, (p)}^{n-1, \, n-1}\in\mbox{Im}\,\bar\partial_t$ for all $t\in\Delta^{\star}$ by (\ref{eqn:half-conclusion3bisbisfinalp}) and by the fact that, thanks to the $\partial\bar\partial$-lemma, $\partial_t\gamma_t^{n-1}\in\mbox{Im}\,\bar\partial_t$ for $t\neq 0$. Recall moreover that the limit of $\bar\partial_t$-exact forms that avoid the {\it small} eigenvalues of $\Delta_t''$ is again $\bar\partial_0$-exact.) Hence $\gamma_0$ is {\it strongly Gauduchon} in this case and the proof of Theorem \ref{The:limitsG} ends here\!; \\

$(2)$\, if  $\partial_t\widetilde\xi_{t, \, (p)}^{n-1, \, n-1}\in\bigoplus\limits_{\lambda\leq\varepsilon_t''}E^{n, n-1}_{\Delta''_t}(\lambda)$ for all $t\in\Delta$, then $\partial_0\widetilde\xi_{0, \, (p)}^{n-1, \, n-1}\in\ker\Delta_0''$ as in the assumption (A1)\!; \\

$(3)$\, if $\partial_t\widetilde\xi_{t, \, (p)}^{n-1, \, n-1}=u_t + v_t$ with $u_t\in\bigoplus\limits_{\lambda\leq\varepsilon_t''}E^{n, n-1}_{\Delta''_t}(\lambda)$ and $v_t\in\bigoplus\limits_{\lambda\geq\varepsilon''}E^{n, n-1}_{\Delta''_t}(\lambda)$ for all $t\in\Delta$, then $u_t, v_t\in\mbox{Im}\,\bar\partial_t$ for all $t\in\Delta^{\star}$, while $u_0\in\ker\Delta_0''$ and $v_0\in\mbox{Im}\, \bar\partial_0$. (In particular $u_0\perp v_0$, hence $||u_0||\leq ||\partial_0\widetilde\xi_{0, \, (p)}^{n-1, \, n-1}||$.) Thus $v_0$ can be absorbed in the $\bar\partial_0$-exact part of $\partial_0\gamma_0^{n-1}$ in (\ref{eqn:half-conclusion3bisbisfinalp}), while the new obstruction $u_0$ to $\partial_0\gamma_0^{n-1}$ being $\bar\partial_0$-exact is $\Delta_0''$-harmonic, much as the former obstruction $\partial_0\widetilde\xi_{0, \, (p)}^{n-1, \, n-1}$ is supposed to be in assumption (A1).  \hfill $\Box$

\vspace{2ex}

 Thus, after possibly replacing $\partial_0\widetilde\xi_{0, \, (p)}^{n-1, \, n-1}$ with $u_0$, we may (and will henceforth) make the assumption (A1). An immediate consequence of (A1) is

\begin{equation}\label{eqn:xi_tilde-nu_delorthog}\ker\Delta_0''\ni\partial_0\widetilde\xi_{0, \, (p)}^{n-1, \, n-1}\perp\partial_0\nu_{0, \, (p)}^{n-1, \, n-1}, \hspace{3ex} \forall\,\,\,\nu_{0, \, (p)}^{n-1, \, n-1}\in (\mbox{Im}\, S_0)\cap{\cal U}_{0, \, (p)},\end{equation}

\noindent because $\partial_0\nu_{0, \, (p)}^{n-1, \, n-1}\in\mbox{Im}\,\bar\partial_0$ by (\ref{eqn:choiceextra}) and because $\ker\Delta_0''\perp\mbox{Im}\,\bar\partial_0$. We get  \\

\noindent $\langle\langle\Delta_0''(\partial_0\widetilde\xi_{0, \, (p)}^{n-1, \, n-1} + \partial_0\nu_{0, \, (p)}^{n-1, \, n-1}), \, \partial_0\widetilde\xi_{0, \, (p)}^{n-1, \, n-1} + \partial_0\nu_{0, \, (p)}^{n-1, \, n-1}\rangle\rangle$ \\
 
\hfill $= \langle\langle\Delta_0''(\partial_0\nu_{0, \, (p)}^{n-1, \, n-1}), \, \partial_0\widetilde\xi_{0, \, (p)}^{n-1, \, n-1}\rangle\rangle + \langle\langle\Delta_0''(\partial_0\nu_{0, \, (p)}^{n-1, \, n-1}), \, \partial_0\nu_{0, \, (p)}^{n-1, \, n-1}\rangle\rangle$ \\

\noindent because $\Delta_0''(\partial_0\widetilde\xi_{0, \, (p)}^{n-1, \, n-1})=0$ by assumption (A1). Now $\Delta_0''(\partial_0\nu_{0, \, (p)}^{n-1, \, n-1})\in\mbox{Im}\,\bar\partial_0$ since $\partial_0\nu_{0, \, (p)}^{n-1, \, n-1}\in\mbox{Im}\,\bar\partial_0$ by the choice (\ref{eqn:choiceextra}) and since $\bar\partial_0$ and $\Delta_0''$ commute. Meanwhile, $\partial_0\widetilde\xi_{0, \, (p)}^{n-1, \, n-1}\in\ker\Delta_0''$ by assumption (A1). Since $\ker\Delta_0''\perp\mbox{Im}\,\bar\partial_0$, the first term on the second line above vanishes. On the other hand, again by the choice (\ref{eqn:choiceextra}) and the definition of ${\cal E}_0^{n, \, n-1}$, we have $\partial_0\nu_{0, \, (p)}^{n-1, \, n-1}\in\bigoplus\limits_{\lambda\geq\varepsilon''}E_{\Delta_0''}^{n, \, n-1}(\lambda)$. It follows that the second term on the second line above satisfies

$$\langle\langle\Delta_0''(\partial_0\nu_{0, \, (p)}^{n-1, \, n-1}), \, \partial_0\nu_{0, \, (p)}^{n-1, \, n-1}\rangle\rangle\geq\varepsilon''\,||\partial_0\nu_{0, \, (p)}^{n-1, \, n-1}||^2,$$

\noindent so we get \\

\noindent $\langle\langle\Delta_0''(\partial_0\widetilde\xi_{0, \, (p)}^{n-1, \, n-1} + \partial_0\nu_{0, \, (p)}^{n-1, \, n-1}), \, \partial_0\widetilde\xi_{0, \, (p)}^{n-1, \, n-1} + \partial_0\nu_{0, \, (p)}^{n-1, \, n-1}\rangle\rangle$

\begin{eqnarray}\label{eqn:delta''zero-lest}\nonumber & \geq & \varepsilon''\, ||\partial_0\nu_{0, \, (p)}^{n-1, \, n-1}||^2 = \varepsilon''\, (||\partial_0\nu_{0, \, (p)}^{n-1, \, n-1}||^2 + ||\partial_0^{\star}\nu_{0, \, (p)}^{n-1, \, n-1}||^2) \\
   & = & \varepsilon''\, \langle\langle\Delta_0'\nu_{0, \, (p)}^{n-1, \, n-1}), \,\nu_{0, \, (p)}^{n-1, \, n-1}\rangle\rangle \geq \varepsilon'\,\varepsilon''\,||\nu_{0, \, (p)}^{n-1, \, n-1}||^2.\end{eqnarray}

\noindent The equality on the second line of (\ref{eqn:delta''zero-lest}) follows from $\partial_0^{\star}\nu_{0, \, (p)}^{n-1, \, n-1}=0$ which in turn follows from $\nu_{0, \, (p)}^{n-1, \, n-1}\in\mbox{Im}\,\partial_0^{\star}\subset\ker\partial_0^{\star}$. (Recall that $\nu_{0, \, (p)}^{n-1, \, n-1}$ has been chosen to have {\it minimal $L^2$-norm} among the $\partial_0$-potentials of $\eta_{0, \, (p)}$ in the definition (\ref{eqn:defS_t}) of the map $S_0$.) The last inequality on the third line in (\ref{eqn:delta''zero-lest}) follows from $\nu_{0, \, (p)}^{n-1, \, n-1}\in\bigoplus\limits_{\mu\geq\varepsilon'}E_{\Delta_0'}^{n-1, \, n-1}(\mu)$ (see (\ref{eqn:u-ball-in}) for $t=0$).

\begin{Concl}\label{Concl:delta''zero-est1} Under the assumption (A1), we have:

\begin{equation}\label{eqn:delta''zero-est1}\langle\langle\Delta_0''(\partial_0\widetilde\xi_{0, \, (p)}^{n-1, \, n-1} + \partial_0\nu_{0, \, (p)}^{n-1, \, n-1}), \, \partial_0\widetilde\xi_{0, \, (p)}^{n-1, \, n-1} + \partial_0\nu_{0, \, (p)}^{n-1, \, n-1}\rangle\rangle\geq\varepsilon'\,\varepsilon''\,||\nu_{0, \, (p)}^{n-1, \, n-1}||^2\end{equation}

\noindent for all $\nu_{0, \, (p)}^{n-1, \, n-1}\in (\mbox{Im}\,S_0)\cap{\cal U}_{0, \, (p)}$.

\end{Concl}

 Now recall that by Conclusion \ref{Concl:sofar} it is only in {\it Case} $1$ that condition $(c)$ of $(\star_p)$ has yet to be obtained. (We have already argued that $(a), (b), (c)$ are simultaneously satisfied in {\it Case} $2$ with the choices  made so far.) Let

$$\alpha_{(p)}:=||\widetilde\xi_{0, \, (p)}|| - ||\widetilde\xi_{0, \, (p)}^{n-1, \, n-1}||>0  \hspace{6ex} \mbox{in \mbox{\it Case} $1$}.$$

\begin{Lem}\label{Lem:max-choice-in-ball} If $\nu_{0, \, (p)}^{n-1, \, n-1}\in (\mbox{Im}\,S_0)\cap{\cal U}_{0, \, (p)}$ is chosen of {\bf maximal $L^2$-norm} among the forms in the intersection of the subspace $\mbox{Im}\,S_0$ with the ball ${\cal U}_{0, \, (p)}$, we have

\begin{equation}\label{eqn:max-choice-in-ball}||\nu_{0, \, (p)}^{n-1, \, n-1}||\geq\alpha_{(p)}\end{equation}

\noindent and $\alpha_{(p)}>0$ in \mbox{\it Case} $1$.

\end{Lem}

\noindent {\it Proof.} It is clear that $\alpha_{(p)}>0$ in {\it Case} $1$ and $\alpha_{(p)}=0$ in {\it Case} $2$.

 In the ball ${\cal U}_{0, \, (p)}$, the ray $R_{(p)}$ emanating from the centre $-\widetilde\xi_{0, \, (p)}^{n-1, \, n-1}$ of ${\cal U}_{0, \, (p)}$ and going through the origin $0\in{\cal U}_{0, \, (p)}$ of the ambient vector space $C^{\infty}_{n-1, \, n-1}(X_0, \, \C)$ cuts the boundary sphere of ${\cal U}_{0, \, (p)}$ in a point that we call $A_{(p)}$. If $d_{(p)}$ denotes the distance from $0$ to $A_{(p)}$, then $d_{(p)}=\alpha_{(p)}$. Meanwhile, the hyperplane $H_{0, \, (p)}$ is orthogonal to the ray $R_{(p)}$ at $0$ and the {\it maximal} $L^2$-norm that a vector $\nu_{0, \, (p)}^{n-1, \, n-1}\in (\mbox{Im}\, S_0)\cap{\cal U}_{0, \, (p)}$ can have attains its {\it minimal} value when $\mbox{Im}\, S_0$ is contained in $H_{0, \, (p)}$. When $\mbox{Im}\, S_0\subset H_{0, \, (p)}$, the vector $\nu_{0, \, (p)}^{n-1, \, n-1}$ can be chosen in the intersection of $\mbox{Im}\, S_0$ with the boundary sphere of ${\cal U}_{0, \, (p)}$ to attain the maximal value that the $L^2$-norm of a vector in $(\mbox{Im}\, S_0)\cap{\cal U}_{0, \, (p)}$ can have in this case. Then in the right-angled triangle formed by the points $0$, $\nu_{0, \, (p)}^{n-1, \, n-1}$ and $A_{(p)}$, the side joining $0$ to $\nu_{0, \, (p)}^{n-1, \, n-1}$ (of length $||\nu_{0, \, (p)}^{n-1, \, n-1}||$) cannot be shorter than the side joining $0$ to $A_{(p)}$ (of length $d_{(p)}=\alpha_{(p)}$) since the angle facing the former side is $\geq \pi/4$ while the angle facing the latter side is $\leq \pi/4$.  \hfill $\Box$

\vspace{2ex}

 Now recall that in {\it Case} $2$ we have $\alpha_{(p)}=0$ and we can choose $\nu_{0, \, (p)}^{n-1, \, n-1}=0$ because $\widetilde\xi_{0, \, (p)}$ already satisfies condition $(c)$ of $(\star_p)$ with $\nu_{0, \, (p)}^{n-1, \, n-1}=0$ for the uniform $\varepsilon_0>0$ obtained at the previous induction steps $1, \dots, p$. Therefore, if in {\it Case} $1$ $\alpha_{(p)}\downarrow 0$ as $p\rightarrow +\infty$, we can satisfy condition $(c)$ of $(\star_p)$ with the uniform $\varepsilon_0>0$ of $(\star_{p_0})(c)$ for all $p\geq p_0$ and for some $p_0\in\N$. \\

 Thus it remains to treat the case covered by the following \\

\noindent {\bf Assumption (A2):} \hspace{6ex} $\alpha_{(p)}\geq\alpha_0>0, \hspace{3ex} \forall p\in\N,$

\noindent {\it for some $\alpha_0>0$ independent of $p\in\N$}.

\vspace{2ex}

 In this case, we get from the estimate (\ref{eqn:delta''zero-est1}) of Conclusion \ref{Concl:delta''zero-est1} and from the estimate (\ref{eqn:max-choice-in-ball}) of Lemma \ref{Lem:max-choice-in-ball} the following

\begin{Concl}\label{Concl:delta''zero-est2} Under the assumptions (A1) and (A2), we have:

\begin{equation}\label{eqn:delta''zero-est2}\langle\langle\Delta_0''(\partial_0\widetilde\xi_{0, \, (p)}^{n-1, \, n-1} + \partial_0\nu_{0, \, (p)}^{n-1, \, n-1}), \, \partial_0\widetilde\xi_{0, \, (p)}^{n-1, \, n-1} + \partial_0\nu_{0, \, (p)}^{n-1, \, n-1}\rangle\rangle\geq\varepsilon'\,\varepsilon''\,\alpha_0^2\end{equation}

\noindent for some $\nu_{0, \, (p)}^{n-1, \, n-1}\in (\mbox{Im}\,S_0)\cap{\cal U}_{0, \, (p)}$ chosen to maximise the $L^2$-norm $||\nu_{0, \, (p)}^{n-1, \, n-1}||$.

\end{Concl}

\vspace{2ex}

 We have thus achieved our purpose of proving the existence of auxiliary forms $\nu_{t, \, (p)}^{n-1, \, n-1}\in (\mbox{Im}\,S_t)\cap{\cal U}_{t, \, (p)}$ (i.e. satisfying (\ref{eqn:ball-inter-ImS}) which automatically guarantees $(a)$ and $(b)$ of $(\star_p)$) such that the numerator of $(c)$ in $(\star_p)$ is uniformly bounded below by a positive constant in {\it Case} $1$.

\vspace{3ex}

\noindent $(II)$\,\, {\it Uniformly bounding the fraction of $(c)$ in $(\star_p)$ from below in Case $1$} \\

  Recall that under $(I)$ above we have been working under the induction hypothesis that the induction steps $1, \dots, p$ had been run and have shown as a result the existence of auxiliary forms $\eta_{t, \, (p)}=\partial_t\nu_{t, \, (p)}^{n-1, \, n-1}=\bar\partial_t\vartheta_{t, \, (p)}^{n, \, n-2}$ satisfying conditions $(a)$, $(b)$ of $(\star_p)$ and (\ref{eqn:delta''zero-est2}). Thus the inductively constructed auxiliary forms satisfy $(a)$ and $(b)$ of $(\star_p)$ for all $p\in\N$ as well as the uniform lower bound:

\begin{equation}\label{eqn:delta''zero-est2-t}\langle\langle\Delta_t''(\partial_t\widetilde\xi_{t, \, (p)}^{n-1, \, n-1} + \partial_t\nu_{t, \, (p)}^{n-1, \, n-1}), \, \partial_t\widetilde\xi_{t, \, (p)}^{n-1, \, n-1} + \partial_t\nu_{t, \, (p)}^{n-1, \, n-1}\rangle\rangle\geq\delta>0,\end{equation}

\noindent for all $t\in\Delta$ (after possibly shrinking $\Delta$ about $0$) and all $p\in\N$, where we have denoted $\delta:=\varepsilon'\,\varepsilon''\,\alpha_0^2>0$ (independent of $t$ and $p$, cf. (\ref{eqn:delta''zero-est2})).

 Now we have:

\begin{eqnarray}\label{eqn:AtpA'tp}\nonumber A_{t, \, (p)}: & = & \langle\langle\Delta_t(\partial_t\widetilde\xi_{t, \, (p)}^{n-1, \, n-1} + \partial_t\nu_{t, \, (p)}^{n-1, \, n-1}), \, \partial_t\widetilde\xi_{t, \, (p)}^{n-1, \, n-1} + \partial_t\nu_{t, \, (p)}^{n-1, \, n-1}\rangle\rangle\\
\nonumber  & = & A'_{t, \, (p)} + \langle\langle\Delta_t''(\partial_t\widetilde\xi_{t, \, (p)}^{n-1, \, n-1} + \partial_t\nu_{t, \, (p)}^{n-1, \, n-1}), \, \partial_t\widetilde\xi_{t, \, (p)}^{n-1, \, n-1} + \partial_t\nu_{t, \, (p)}^{n-1, \, n-1}\rangle\rangle\\
   & \geq & A'_{t, \, (p)} + \delta, \hspace{3ex} t\in\Delta, \,\, p\in\N,\end{eqnarray}

\noindent where we have denoted $A'_{t, \, (p)}\!\!:= \langle\langle\Delta'_t(\partial_t\widetilde\xi_{t, \, (p)}^{n-1, \, n-1} + \partial_t\nu_{t, \, (p)}^{n-1, \, n-1}), \, \partial_t\widetilde\xi_{t, \, (p)}^{n-1, \, n-1} + \partial_t\nu_{t, \, (p)}^{n-1, \, n-1}\rangle\rangle$ and by $A_{t, \, (p)}$ the analogous expression with $\Delta_t$ in place of $\Delta_t'$. (To justify the identity between the top two lines in (\ref{eqn:AtpA'tp}), recall that for any {\it pure-type} form $u$ one has $\langle\langle\Delta_tu, \, u\rangle\rangle = \langle\langle\Delta_t'u, \, u\rangle\rangle + \langle\langle\Delta_t''u, \, u\rangle\rangle$ by (\ref{eqn:proof-laplace-comparison}).)

 Recall that in order to guarantee condition $(c)$ of $(\star_p)$ for all $p\in\N$ we need to prove the existence of an $\varepsilon_0>0$ independent of both $t\in\Delta$ and $p\in\N$ such that

\begin{equation}\label{eqn:need-end-c-AA'}A_{t, \, (p)}\geq (1+\varepsilon_0)\, A'_{t, \, (p)}, \hspace{3ex} t\in\Delta, \,\, p\in\N.\end{equation}

\noindent Since (\ref{eqn:AtpA'tp}) holds, it suffices to get a uniform $\varepsilon_0>0$ as above such that

\begin{equation}\label{need-end-c-frac}A'_{t, \, (p)} + \delta \geq (1+\varepsilon_0)\, A'_{t, \, (p)} \hspace{3ex} \mbox{or equivalently} \hspace{3ex} A'_{t, \, (p)}\leq\frac{\delta}{\varepsilon_0}, \hspace{3ex} t\in\Delta, \,\, p\in\N.\end{equation}

\noindent  The existence of such a uniform $\varepsilon_0>0$ is of course guaranteed if we can prove that $A'_{t, \, (p)}$ is uniformly bounded above. Since $A'_{t, \, (p)}\leq A_{t, \, (p)}$, it suffices to prove the existence of a uniform upper bound for the latter quantity.

\begin{Lem}\label{Lem:need-end-c-upperb} In the above notation, the auxiliary forms $\eta_{t, \, (p)}=\partial_t\nu_{t, \, (p)}^{n-1, \, n-1}=\bar\partial_t\vartheta_{t, \, (p)}^{n, \, n-2}\in{\cal E}_t^{n, \, n-1}$ constructed by the induction procedure set up in the preceding paragraphs and with the choices made there satisfy

\begin{equation}\label{eqn:need-end-c-upperb}A'_{t, \, (p)}\leq A_{t, \, (p)}\leq M < +\infty, \hspace{3ex} t\in\Delta, \,\, p\in\N,\end{equation}

\noindent for some $M$ independent of both $t\in\Delta$ and $p\in\N$.

\end{Lem}

\noindent {\it Proof.} Recall that in the induction process we solve the equations (cf. (\ref{eqn:re-iterative-pot})):

\begin{equation}\label{eqn:ind-eqn}d\widetilde\xi_{t, \, (p+1)}=\partial_t(\widetilde\xi_{t, \, (p)}^{n-1, \, n-1}+\nu_{t, \, (p)}^{n-1, \, n-1}), \hspace{3ex} t\in\Delta, \,\, p\in\N,\end{equation}

\noindent and we choose $\widetilde\xi_{t, \, (p+1)}$ to be the minimal $L^2$-norm solution for every given $p\in\N$. Thus

\begin{equation}\label{eqn:ind-minsol}\widetilde\xi_{t, \, (p+1)}=\Delta_t^{-1}d_t^{\star}\partial_t(\widetilde\xi_{t, \, (p)}^{n-1, \, n-1}+\nu_{t, \, (p)}^{n-1, \, n-1}), \hspace{3ex} t\in\Delta, \,\, p\in\N,\end{equation}

\noindent and

\begin{eqnarray}\label{eqn:Btp}\nonumber||\widetilde\xi_{t, \, (p+1)}||^2 & = & ||\Delta_t^{-\frac{1}{2}}\partial_t(\widetilde\xi_{t, \, (p)}^{n-1, \, n-1}+\nu_{t, \, (p)}^{n-1, \, n-1})||^2\\
  & = & B_{t, \, (p)}, \hspace{3ex} t\in\Delta, \,\, p\in\N,\end{eqnarray}

\noindent where we have denoted

$$B_{t, \, (p)}:=\langle\langle\Delta_t^{-1}\partial_t(\widetilde\xi_{t, \, (p)}^{n-1, \, n-1} + \nu_{t, \, (p)}^{n-1, \, n-1}), \, \partial_t(\widetilde\xi_{t, \, (p)}^{n-1, \, n-1} + \nu_{t, \, (p)}^{n-1, \, n-1})\rangle\rangle.$$

\noindent It is clear that if $B_{t, \, (p)}$ became arbitrarily small when $p\rightarrow +\infty$, then $||\widetilde\xi_{t, \, (p+1)}||$ would become arbitrarily small. This would give right away the conclusion of Corollary \ref{Cor:finalnormest} below and the proof of Theorem \ref{The:limitsG} would follow as explained at the end of the paper. This gives a hint that $A_{t, \, (p)}$ is likely to satisfy the uniform upper bound (\ref{eqn:need-end-c-upperb}) at least in the complementary case (i.e. when $B_{t, \, (p)}$ is uniformly bounded below by a positive constant). Here are the details.

 If we denote $\varpi_{t, \, (p)}:=\partial_t(\widetilde\xi_{t, \, (p)}^{n-1, \, n-1} + \nu_{t, \, (p)}^{n-1, \, n-1})$, we know that

$$\varpi_{t, \, (p)}=d\,\widetilde\xi_{t, \, (p+1)}, \hspace{3ex} \mbox{with}\hspace{2ex} \widetilde\xi_{t, \, (p+1)}\in\mbox{Im}\, d_t^{\star}\subset\ker\,d_t^{\star},\,\, t\in\Delta,\,\, p\in\N.$$

\noindent So we get

\begin{eqnarray}\label{eqn:Atp-Laplacian}\nonumber A_{t, \, (p)} & = & \langle\langle\Delta_t\varpi_{t, \, (p)}, \, \varpi_{t, \, (p)}\rangle\rangle \\
 \nonumber & = & ||d\,\varpi_{t, \, (p)}||^2 + ||d_t^{\star}\varpi_{t, \, (p)}||^2 = ||d_t^{\star}\varpi_{t, \, (p)}||^2 \\
\nonumber  & = & ||d_t^{\star}d\,\widetilde\xi_{t, \, (p+1)}||^2 = ||d_t^{\star}d\,\widetilde\xi_{t, \, (p+1)} + d\,d_t^{\star}\widetilde\xi_{t, \, (p+1)}||^2\\
  & = & ||\Delta_t\,\widetilde\xi_{t, \, (p+1)}||^2, \hspace{3ex} t\in\Delta,\,\, p\in\N.\end{eqnarray}

 Now observe that the proof of Lemma \ref{Lem:decrease} shows that the families of forms $(\widetilde\xi_{t, \, (p)})_{t\in\Delta}$ ($p\in\N$) defined by solving equations (\ref{eqn:re-iterative-pot}) for $p-1$ satisfy inequalities analogous to the inequalities (\ref{eqn:decrease}) for $(\xi_{t, \, (p)})_{t\in\Delta}$ ($p\in\N$):

\begin{equation}\label{eqn:xi-Omega-ref-re}||\widetilde\xi_{t, \, (p+1)}||\leq ||\widetilde\Omega_{t, \, (p)}^{n-1, \, n-1}||, \hspace{3ex} t\in\Delta,\,\, p\in\N, \hspace{3ex} (\mbox{cf.}\,\, (\ref{eqn:xi-Omega-ref}))\end{equation}

\noindent by comparison of the minimal $d$ and $\partial_t$-potentials of the given form $\varpi_{t, \, (p)}$,

\begin{equation}\label{eqn:L2norm-min-del-re}||\widetilde\Omega_{t, \, (p)}^{n-1, \, n-1}||\leq ||\widetilde\xi_{t, \, (p)}^{n-1, \, n-1} + \nu_{t, \, (p)}^{n-1, \, n-1}||, \hspace{3ex} t\in\Delta,\,\, p\in\N, \hspace{3ex} (\mbox{cf.}\,\, (\ref{eqn:L2norm-min-del}))\end{equation}

\noindent by minimality of $\widetilde\Omega_{t, \, (p)}^{n-1, \, n-1}$ among the $\partial_t$-potentials of $\varpi_{t, \, (p)}$, and

\begin{equation}\label{eqn:bstarp-re}||\widetilde\xi_{t, \, (p)}^{n-1, \, n-1} + \nu_{t, \, (p)}^{n-1, \, n-1}||\leq ||\widetilde\xi_{t, \, (p)}||, \hspace{3ex} t\in\Delta,\,\, p\in\N,\end{equation}

\noindent by $(b)$ of $(\star_p)$. The last three inequalities add up to

\begin{equation}\label{eqn:tilde-decrease}||\widetilde\xi_{t, \, (p+1)}||\leq ||\widetilde\xi_{t, \, (p)}||, \hspace{3ex} t\in\Delta,\,\, p\in\N, \hspace{3ex} (\mbox{cf.}\,\, (\ref{eqn:decrease})).\end{equation}

\noindent The sequence $(||\widetilde\xi_{t, \, (p)}||)_{p\in\N}$ is thus non-increasing (hence bounded above) for each $t\in\Delta$. After slightly shrinking $\Delta$ about $0$, let

\begin{equation}\label{eqn:xi-tilde-upperb}M_1\!\!:=\sup\limits_{t\in\Delta, \, p\in\N}||\widetilde\xi_{t, \, (p)}|| = \sup\limits_{t\in\Delta}||\widetilde\xi_{t, \, (0)}|| = \sup\limits_{t\in\Delta}||\xi_t|| <+\infty.\end{equation}

\vspace{2ex}

 Recall that in view of formula (\ref{eqn:Atp-Laplacian}) we need to show that $||\Delta_t\widetilde\xi_{t, \, (p)}||$ is bounded above independently of $t\in\Delta$ and $p\in\N$. Only the uniform boundedness w.r.t. $p$ has yet to be justified. Note that $\Delta_t$ does not depend on $p$. We need a slight refinement of (\ref{eqn:tilde-decrease}).

 For every $t\in\Delta$ and $p\in\N$ let

\begin{equation}\label{eqn:xi-tilde-delta-decomp}\widetilde\xi_{t, \, (p)}=\sum\limits_{j=0}^{+\infty}u_j^{(p)}(t), \hspace{3ex} \mbox{with}\hspace{2ex} u_j^{(p)}(t)\in E_{\Delta_t}(\lambda_j),\end{equation}

\noindent be the decomposition of $\widetilde\xi_{t, \, (p)}$ w.r.t. the eigenspaces $E_{\Delta_t}(\lambda_j)$ of $\Delta_t$. The eigenvalues $\lambda_j=\lambda_j(t)$ of $\Delta_t$, ordered (without repetitions) increasingly, tend to $+\infty$ as $j$ tends to $+\infty$. Inequality (\ref{eqn:tilde-decrease}) translates to

\begin{equation}\label{eqn:xi-tilde-decomp-ineq}||\widetilde\xi_{t, \, (p+1)}||^2=\sum\limits_{j=0}^{+\infty}||u_j^{(p+1)}(t)||^2\leq\sum\limits_{j=0}^{+\infty}||u_j^{(p)}(t)||^2=||\widetilde\xi_{t, \, (p)}||^2, \hspace{3ex} t\in\Delta,\,\, p\in\N.\end{equation}

\noindent Meanwhile, we clearly have

\begin{equation}\label{eqn:xi-tilde-delta-norm}||\Delta_t\widetilde\xi_{t, \, (p)}||^2=\sum\limits_{j=0}^{+\infty}\lambda_j^2\,||u_j^{(p)}(t)||^2, \hspace{3ex} t\in\Delta,\,\, p\in\N.\end{equation}

 The inductive process that produced the forms $(\widetilde\xi_{t, \, (p)})$ shows in effect that the norm inequality (\ref{eqn:tilde-decrease}) occurs {\it component-wise}, i.e. for every $j\in\N$ we have:

\begin{equation}\label{eqn:c-wise-xi-tilde}||u_j^{(p+1)}(t)||\leq ||u_j^{(p)}(t)||,  \hspace{3ex} t\in\Delta,\,\, p\in\N.\end{equation}

\noindent Indeed, recall that by (\ref{eqn:laplace-comparison}) any pure-type form $u$ satisfies $\langle\langle\Delta_tu, \, u\rangle\rangle\geq\langle\langle\Delta_t'u, \, u\rangle\rangle$, hence inequality (\ref{eqn:xi-Omega-ref-re}) occurs {\it component-wise}. Inequality (\ref{eqn:L2norm-min-del-re}) occurs {\it component-wise} as well since $\widetilde\xi_{t, \, (p)}^{n-1, \, n-1} + \nu_{t, \, (p)}^{n-1, \, n-1}$ is obtained from $\widetilde\Omega_{t, \, (p)}^{n-1, \, n-1}$ by adding a form (lying in $\ker\partial_t$) that is orthogonal to the minimal $L^2$-norm $\partial_t$-potential $\widetilde\Omega_{t, \, (p)}^{n-1, \, n-1}\in\mbox{Im}\,\partial_t^{\star}\perp\ker\partial_t$. On the other hand, $\nu_{t, \, (p)}^{n-1, \, n-1}$ is chosen to lie in $\mbox{Im}\,S_t\cap{\cal U}_{t, \, (p)}$ (by (\ref{eqn:ball-inter-ImS})) and to have {\it maximal} $L^2$-norm among these forms (by Lemma \ref{Lem:max-choice-in-ball}) while $S_t$ is independent of $p$ and the radius of the ball ${\cal U}_{t, \, (p)}$ is non-increasing w.r.t. $p\in\N$ by (\ref{eqn:tilde-decrease}). Hence we can choose the forms $\nu_{t, \, (p)}^{n-1, \, n-1}$ such that 

$$||\nu_{t, \, (p)}^{n-1, \, n-1}||\leq ||\nu_{t, \, (p-1)}^{n-1, \, n-1}  || \hspace{6ex} \mbox{\it component-wise}, \hspace{3ex} t\in\Delta,\,\, p\in\N^{\star}.$$

\noindent Thus we obtain (\ref{eqn:c-wise-xi-tilde}) inductively on $p\in\N$: if (\ref{eqn:c-wise-xi-tilde}) has been shown for $p-1$, then for all $t\in\Delta$  we have

$$||\widetilde\xi_{t, \, (p)}^{n-1, \, n-1} + \nu_{t, \, (p)}^{n-1, \, n-1}||\leq ||\widetilde\xi_{t, \, (p-1)}^{n-1, \, n-1} + \nu_{t, \, (p-1)}^{n-1, \, n-1}  || \hspace{6ex} \mbox{\it component-wise},$$

\noindent which implies $||\widetilde\xi_{t, \, (p+1)}||\leq ||\widetilde\xi_{t, \, (p)}||$ {\it component-wise}. This is nothing but (\ref{eqn:c-wise-xi-tilde}).

 Now (\ref{eqn:xi-tilde-delta-norm}) and (\ref{eqn:c-wise-xi-tilde}) combine to show the existence of a uniform upper bound for the Laplacian of $\widetilde\xi_{t, \, (p)}$ (after slightly shrinking $\Delta$ about $0$):

\begin{equation}\label{eqn:xi-tilde-upperb}M\!\!:=\sup\limits_{t\in\Delta, \, p\in\N}||\Delta_t\,\widetilde\xi_{t, \, (p+1)}||<+\infty,\end{equation}

\noindent which in view of (\ref{eqn:Atp-Laplacian}) is nothing but (\ref{eqn:need-end-c-upperb}).

 Lemma \ref{Lem:need-end-c-upperb} is proved.  \hfill $\Box$

\vspace{3ex}

 We can now explicitly achieve (\ref{need-end-c-frac}), hence also (\ref{eqn:need-end-c-AA'}) which is equivalent to condition $(c)$ of $(\star_p)$. Indeed, estimate (\ref{eqn:need-end-c-upperb}) obtained in Lemma \ref{Lem:need-end-c-upperb} shows that the inductively constructed auxiliary forms fulfill condition $(c)$ of $(\star_p)$ with the uniform $\varepsilon_0>0$ defined by

$$\varepsilon_0:=\frac{\delta}{M} = \frac{\varepsilon'\,\varepsilon''\,\alpha_0^2}{M}>0,$$

\noindent where $\delta:=\varepsilon'\,\varepsilon''\,\alpha_0^2>0$ is the uniform lower bound of (\ref{eqn:delta''zero-est2-t}) and $M<+\infty$ is the uniform upper bound of (\ref{eqn:xi-tilde-upperb}).

\vspace{2ex}

The existence of the auxiliary forms is thus accounted for.   \hfill $\Box$

\section {Final arguments in proving Theorem \ref{The:limitsG}}\label{section:end-proof}

 With the new inductive construction based on auxiliary forms in place, the identities of Lemma \ref{Lem:p-iterations} obeyed by $\xi_{t, \, (p)}$ are transformed into the following identities obeyed by $\widetilde\xi_{t, \, (p)}$.

\begin{Lem}\label{Lem:re-p-iterations} The family $(\widetilde{\xi}_{t, \, (p)})_{t\in\Delta}$ of $(2n-2)$-forms constructed above varies in a $C^{\infty}$ way with $t$ (up to $t=0$) and satisfies for all $t\in\Delta$ and all $p\in\N$:

\begin{eqnarray}\label{eqn:re-stepp}\nonumber\partial_t(\gamma_t^{n-1} & - & \widetilde{\xi}_{t, \, (p)}^{n-1, \, n-1} - \overline{\widetilde{\xi}_{t, \, (p)}^{n-1, \, n-1}}) = \bar\partial_t(\widetilde{\xi}_{t, \, (p)}^{n, \, n-2} + \overline{\widetilde{\xi}_{t, \, (p)}^{n-2, \, n}} + \dots + \widetilde{\xi}_{t, \, (1)}^{n, \, n-2} + \xi_t^{n, \, n-2} \\
 & + & w_t - \vartheta_t^{n, \, n-2} - \vartheta_{t, \, (1)}^{n, \, n-2} - \dots - \vartheta_{t, \, (p-1)}^{n, \, n-2}).\end{eqnarray}

\end{Lem}

\noindent {\it Proof.} It follows trivially from (\ref{eqn:half-conclusion3bisbisfinal}) with $p+1$ replaced by $p$ and the fact that $d\, \widetilde{\xi}_{t, \, (p)}=\partial_t\widetilde{\xi}_{t, \, (p)}^{n-1, \, n-1} + \bar\partial_t\widetilde{\xi}_{t, \, (p)}^{n, \, n-2}$ is of type $(n, \, n-1)$ (thus its $(n-1, \, n)$-component vanishes, hence $-\bar\partial_t\widetilde{\xi}_{t, \, (p)}^{n-1, \, n-1} = \partial_t\widetilde{\xi}_{t, \, (p)}^{n-2, \, n}$ and taking conjugates $-\partial_t\overline{\widetilde{\xi}_{t, \, (p)}^{n-1, \, n-1}} = \bar\partial_t\overline{\widetilde{\xi}_{t, \, (p)}^{n-2, \, n}}$) by arguments analogous to those of Lemma \ref{Lem:p-iterations}. \hfill $\Box$

\vspace{2ex}

 The next, more substantial step is to show that the $L^2$-norm of $\widetilde{\xi}_{t, \, (p)}^{n-1, \, n-1}$ decreases strictly at each step $p$ of the above inductive construction in a way that guarantees it to become arbitrarily small when $p$ becomes large enough. The following lemma and its corollary provide the final argument to the proof of Theorem \ref{The:limitsG} and, implicitly, to that of Theorem \ref{The:projdef}.

\begin{Lem}\label{Lem:finalnormest} There exists $\varepsilon>0$ independent of $t\in\Delta$ and of $p\in\N$ such that the minimal $L^2$-norm solutions $\widetilde{\Omega}_{t, \, (p)}^{n-1, \, n-1}$ and $\widetilde{\xi}_{t, \, (p+1)}$ of the equations 

\begin{equation}\label{eqn:iterative-pot-re}\partial_t\widetilde{\Omega}_{t, \, (p)}^{n-1,\, n-1} = \partial_t(\widetilde{\xi}_{t, \, (p)}^{n-1, \, n-1} + \nu_{t, \, (p)}^{n-1, \, n-1}) \hspace{2ex} \mbox{and} \hspace{2ex} d\, \widetilde{\xi}_{t,\, (p+1)}= \partial_t(\widetilde{\xi}_{t, \, (p)}^{n-1, \, n-1} + \nu_{t, \, (p)}^{n-1, \, n-1})\end{equation}

\noindent satisfy the $L^2$-norm estimates:

\begin{equation}\label{eqn:sol-comp}||\widetilde{\xi}_{t, \, (p+1)}||\leq \frac{1}{\sqrt{1+\varepsilon}}\, ||\widetilde{\Omega}_{t, \, (p)}^{n-1, \, n-1}||, \hspace{2ex} t\in\Delta, \, p\in\N.\end{equation}

\end{Lem}

 Before proving this statement, we notice an immediate corollary.

\begin{Cor}\label{Cor:finalnormest} The forms $\widetilde{\xi}_{t, \, (p)}$ obtained above satisfy

\begin{equation}\label{eqn:iteration-compp}||\widetilde{\xi}_{t, \, (p)}||\leq \frac{1}{(\sqrt{1+\varepsilon})^p}\, ||\xi_t||, \hspace{2ex} t\in\Delta, \,\, p\in\N.\end{equation}

 In particular, $||\widetilde{\xi}_{t, \, (p)}||$ (hence also $||\widetilde{\xi}_{t, \, (p)}^{n-1, \, n-1}||$ which is $\leq ||\widetilde{\xi}_{t, \, (p)}||$) becomes arbitrarily small, uniformly w.r.t. $t\in\Delta$ and $p\gg 1$, if the number $p\in\N$ of iterations is sufficiently large.

\end{Cor}

\noindent {\it Proof of Corollary \ref{Cor:finalnormest}.} From Lemma \ref{Lem:finalnormest} we get the following inequalities:

$$||\widetilde{\xi}_{t, \, (p+1)}|| \leq \frac{1}{\sqrt{1+\varepsilon}}\, ||\widetilde{\Omega}_{t, \, (p)}^{n-1, \, n-1}||\leq \frac{1}{\sqrt{1+\varepsilon}}\, ||\widetilde{\xi}_{t, (p)}^{n-1, \, n-1} + \nu_{t, \, (p)}^{n-1, \, n-1}||, \hspace{2ex} p\in\N.$$

The latter inequality follows from the $L^2$-norm minimality of $\widetilde{\Omega}_{t, \, (p)} ^{n-1, \, n-1}$ among the solutions of the equation $\partial_t\widetilde{\Omega}_{t, \, (p)}^{n-1, \, n-1}=\partial_t(\widetilde{\xi}_{t, \, (p)}^{n-1, \, n-1}+\nu_{t, \, (p)}^{n-1, \, n-1})$. Combining with $(b)$ of properties $(\star_p)$, we get

$$||\widetilde{\xi}_{t, \, (p+1)}|| \leq \frac{1}{\sqrt{1+\varepsilon}}\, ||\widetilde{\xi}_{t, \, (p)}||, \hspace{2ex} t\in\Delta, \,\, p\in\N.$$

\noindent Letting $p$ run through $0,\dots , p-1$, these inequalities multiply up to (\ref{eqn:iteration-compp}).  \hfill $\Box$

\vspace{2ex}

 We now come to the key task of proving Lemma \ref{Lem:finalnormest}. However, the ground has been largely prepared by Lemma \ref{Lem:decrease} and Observation \ref{Obs:t-unif} whose proofs outlined the difficulties and explained how to overcome them under certain hypotheses, as well as by the construction of auxiliary forms $\eta_{t, \, (p)}$ satisfying conditions $(\star_p)$ which enable those hypotheses to be met. The remaining arguments are almost purely formal.

\vspace{2ex}

\noindent {\it Proof of Lemma \ref{Lem:finalnormest}.} Recall the notation $\widetilde{\xi}_{t, \, (0)}^{n-1, \, n-1}: = \xi_t^{n-1, \, n-1}$, $\widetilde{\Omega}_{t, \, (0)}^{n-1, \, n-1}:= \widetilde{\Omega}_t^{n-1, \, n-1}$ and $\nu_{t, \, (0)}^{n-1, \, n-1}:= \nu_t^{n-1, \, n-1}$. Set $\varpi_{t, \, (p)} := \partial_t(\widetilde{\xi}_{t, \, (p)}^{n-1, \, n-1} + \nu_{t, \, (p)}^{n-1, \, n-1})$, the right-hand term of equations (\ref{eqn:iterative-pot-re}). The minimal $L^2$-norm solutions of equations (\ref{eqn:iterative-pot-re}) are explicitly given by the formulae:

\begin{equation}\label{eqn:minsol-formulae} \widetilde{\Omega}_{t, \, (p)}^{n-1, \, n-1}=\Delta_t^{'-1}\partial_t^{\star}\varpi_{t, \, (p)} \hspace{2ex} \mbox{and} \hspace{2ex} \widetilde{\xi}_{t, \, (p+1)}=\Delta_t^{-1}d_t^{\star}\varpi_{t, \, (p)}, \hspace{3ex} t\in\Delta, \hspace{1ex} p\in\N.\end{equation}

\noindent Thus by (\ref{eqn:delta'-halfnorm}) and (\ref{eqn:delta-halfnorm}) with $u=\varpi_{t, \, (p)}$, the proof of Lemma \ref{Lem:finalnormest} reduces to proving that, for some $\varepsilon>0$ independent of $t\in\Delta$ and $p\in\N$, we have:

\begin{equation}\label{eqn:reduced-finalnormest}||\Delta_t^{-\frac{1}{2}}\varpi_{t, \, (p)}|| \leq \frac{1}{\sqrt{1+\varepsilon}}\, ||\Delta_t^{'-\frac{1}{2}}\varpi_{t, \, (p)}||, \hspace{3ex} t\in\Delta, \,\, p\in\N.\end{equation}

 Now the forms $\eta_{t, \, (p)}=\partial_t\nu_{t, \, (p)}^{n-1, \, n-1}$ have been chosen to satisfy conditions $(\star_p)$ whose part $(c)$ translates to:

\begin{equation}\label{eqn:epsilon-cond}0<\varepsilon_0 \leq \frac{\langle\langle\Delta_t'' \varpi_{t, \, (p)}, \, \varpi_{t, \, (p)}\rangle\rangle}{\langle\langle\Delta_t' \varpi_{t, \, (p)}, \, \varpi_{t, \, (p)}\rangle\rangle},  \hspace{3ex}  t\in\Delta, \,\, p\in\N,\end{equation}

\noindent for an $\varepsilon_0 >0$ independent of both $t\in\Delta$ and $p\in\N$.
 
 By the choice $(a)$ of $(\star_p)$, we have $\eta_{t, \, (p)}=\partial_t\nu_{t, \, (p)}^{n-1, \, n-1}=\bar\partial_t\vartheta_{t, \, (p)}^{n, \, n-2}$, hence $\eta_{t, \, (p)}$ is $\bar\partial_t$-exact for all $t\in\Delta$ and all $p\in\N$. It follows that:\\

 $(i)$\,  the form $\varpi_{t, \, (p)}=\partial_t\widetilde{\xi}_{t, \, (p)}^{n-1, \, n-1}+\eta_{t, \, (p)}$ is $\bar\partial_t$-exact for all $t\neq 0$, hence $\varpi_{t, \, (p)}$ is orthogonal to $\ker\Delta_t''$ for all $t\neq 0$;

$(ii)$\, when $t=0$, the form $\varpi_{0, \, (p)}=\partial_0\widetilde{\xi}_{0, \, (p)}^{n-1, \, n-1}+\eta_{0, \, (p)}$ cannot be $\Delta_0''$-harmonic. 

 Indeed, otherwise the condition $(c)$ of $(\star_p)$ would be violated (see (\ref{eqn:epsilon-cond}) for $t=0$) unless we also have $\Delta_0'\varpi_{0, \, (p)}=0$. However, in this latter case the $\partial_0$-exact form $\varpi_{0, \, (p)}=\partial_0(\widetilde{\xi}_{0, \, (p)}^{n-1, \, n-1} + \nu_{0, \, (p)}^{n-1, \, n-1})$ would have to vanish (since $\mbox{Im}\,\partial_0\perp\ker\Delta_0'$) and $\partial_0\gamma_0^{n-1}$ would be $\bar\partial_0$-exact by (\ref{eqn:half-conclusion3bisbisp+1}) applied at $t=0$. Then $\gamma_0$ would be a {\it strongly Gauduchon} metric on $X_0$ and the proof of Theorem \ref{The:limitsG} would be complete.

 We conclude that, for every fixed $p\in\N$, the family $(\varpi_{t, \, (p)})_{t\in\Delta}$ satisfies the non-$\Delta_t''$-harmonicity hypothesis (\ref{eqn:unif-hyp}), hence also estimate (\ref{eqn:oneplusepsilon-t-est}) uniformly w.r.t. $t\in\Delta$. 

 Moreover, by (\ref{eqn:proof-laplace-comparison}) and by $\varpi_{t, \, (p)}$ being of pure type, we have 

$$\langle\langle\Delta_t \varpi_{t, \, (p)}, \, \varpi_{t, \, (p)}\rangle\rangle = \langle\langle\Delta_t' \varpi_{t, \, (p)}, \, \varpi_{t, \, (p)}\rangle\rangle + \langle\langle\Delta_t'' \varpi_{t, \, (p)}, \, \varpi_{t, \, (p)}\rangle\rangle,$$

\noindent so the uniform estimate (\ref{eqn:epsilon-cond}) amounts to

\begin{equation}\label{eqn:epsilont-pit}\langle\langle\Delta_t \varpi_{t, \, (p)}, \, \varpi_{t, \, (p)}\rangle\rangle \geq (1+\varepsilon_0)\, \langle\langle\Delta_t' \varpi_{t, \, (p)}, \, \varpi_{t, \, (p)}\rangle\rangle, \hspace{2ex} t\in\Delta, \hspace{1ex} p\in\N,\end{equation}

\noindent which provides unifomity w.r.t. $p\in\N$ besides the uniformity w.r.t. $t\in\Delta$. This proves Lemma \ref{Lem:finalnormest}. \hfill $\Box$

\vspace{2ex}

\noindent {\it End of proof of Theorem \ref{The:limitsG}.} By Corollary \ref{Cor:finalnormest}, the $L^2$-norm $||\widetilde{\xi}_{t, \, (p)}^{n-1, \, n-1}||$ can be made arbitrarily small, uniformly with respect to $t\in\Delta$ and $p\gg 1$, if $p$ is chosen sufficiently large. In particular, so can the $L^2$-norm $||\widetilde{\xi}_{0, \, (p)}^{n-1, \, n-1}||$. 

 Thanks to Lemma \ref{Lem:small-correction}, if $p$ is sufficiently large, we get a $C^{\infty}$ positive definite $J_0\!-\!(1, \, 1)$-form $\rho_0>0$ such that 

$$\partial_0\rho_0^{n-1} - \partial_0\bigg(\gamma_0^{n-1}-\widetilde{\xi}_{0, \, (p)}^{n-1, \, n-1}-\overline{\widetilde{\xi}_{0, \, (p)}^{n-1, \, n-1}}\bigg) \in \mbox{Im}\,(\partial_0\bar\partial_0).$$

\noindent Since $\partial_0(\gamma_0^{n-1}-\widetilde{\xi}_{0, \, (p)}^{n-1, \, n-1}-\overline{\widetilde{\xi}_{0, \, (p)}^{n-1, \, n-1}})$ is known to be $\bar\partial_0$-exact by identity (\ref{eqn:re-stepp}) of Lemma \ref{Lem:re-p-iterations}, we see that $\partial_0\rho_0^{n-1}$ must be $\bar\partial_0$-exact, hence $\rho_0$ is a {\it strongly Gauduchon} metric on $X_0$. The proof of Theorem \ref{The:limitsG} is complete. \hfill $\Box$   

\vspace{2ex}

 As explained earlier, Theorem \ref{The:limitsG} combined with Theorem \ref{The:projdef_sG-assumption} proved in [Pop13] proves Theorem \ref{The:projdef}.

\vspace{6ex}

\noindent {\bf References.} \\

\noindent [AK13]\, D. Angella, H. Kasuya --- {\it Cohomologies of Deformations of Solvmanifolds and Closedness of Some Properties} --- arXiv:1305.6709 [math.CV].

\vspace{1ex}

\noindent [CE53]\, E. Calabi, B. Eckmann --- {\it A Class of Compact, Complex Manifolds Which Are Not Algebraic} --- Ann. of Math. {\bf 58} (1953) 494-500.

\vspace{1ex}

\noindent [COUV16]\, M. Ceballos, A. Otal, L. Ugarte, R. Villacampa --- {\it In variant Complex Structures on $6$-Nilmanifolds: Classification, Fr\"olicher Spectral Sequence and Special Hermitian Metrics} --- J. Geom. Anal. {\bf 26} (2016), 252-286.

\vspace{1ex}

\noindent [DGMS75]\, P. Deligne, Ph. Griffiths, J. Morgan, D. Sullivan --- {\it Real Homotopy Theory of K\"ahler Manifolds} --- Invent. Math. {\bf 29} (1975), 245-274.

\vspace{1ex}

\noindent [Ehr47]\, C. Ehresmann ---{\it Sur les espaces fibr\'es diff\'erentiables} -- C. R. Acad. Sci. Paris {\bf 224} (1947), 1611-1612.

\vspace{1ex}

\noindent [FOU15]\, A. Fino, A. Otal, L. Ugarte --- {\it Six-dimensional Solvmanifolds with Holomorphically Trivial Canonical Bundle} --- Int. Math. Res. Not. IMRN 2015, no. 24, 13757-13799.

\vspace{1ex}

\noindent [Fri17]\, R. Friedman --- {\it The $\partial\bar\partial$-Lemma for General Clemens Manifolds} --- arXiv e-print AG 1708.00828v1.

\vspace{1ex}

\noindent [Gau77]\, P. Gauduchon --- {\it Le th\'eor\`eme de l'excentricit\'e nulle} --- C.R. Acad. Sc. Paris, S\'erie A, t. {\bf 285} (1977), 387-390.

\vspace{1ex}

\noindent [Hir62]\, H. Hironaka --- {\it An Example of a Non-K\"ahlerian Complex-Analytic Deformation of K\"ahlerian Complex Structures} --- Ann. of Math. (2) {\bf 75 (1)} (1962), 190-208.

\vspace{1ex}

\noindent [Kod86]\, K. Kodaira --- {\it Complex Manifolds and Deformations of Complex Structures} --- Grundlehren der Math. Wiss. {\bf 283}, Springer (1986).

\vspace{1ex}

\noindent [KS60]\, K. Kodaira, D.C. Spencer --- {\it On Deformations of Complex Analytic Structures, III. Stability Theorems for Complex Structures} --- Ann. of Math. {\bf 71}, no.1 (1960), 43-76.

\vspace{1ex}

\noindent [Mic82]\, M. L. Michelsohn --- {\it On the Existence of Special Metrics in Complex Geometry} --- Acta Math. {\bf 149} (1982), no. 3-4, 261-295.

\vspace{1ex}

\noindent [Pop10]\, D. Popovici --- {\it Limits of Moishezon Manifolds Under Holomorphic Deformations} -- arXiv e-print AG 1003.3605v1.

\vspace{1ex}

\noindent [Pop13]\, D. Popovici --- {\it Deformation Limits of Projective Manifolds: Hodge Numbers and Strongly Gauduchon Metrics} --- Invent. Math. {\bf 194} (2013), 515-534.

\vspace{1ex}

\noindent [Pop14]\,  D. Popovici --- {\it Deformation Openness and Closedness of Various Classes of Compact Complex Manifolds; Examples} --- Ann. Sc. Norm. Super. Pisa Cl. Sci. (5), Vol. XIII (2014), 255-305.

\vspace{1ex}

\noindent [Sch07]\, M. Schweitzer --- {\it Autour de la cohomologie de Bott-Chern} --- arXiv e-print math.AG/0709.3528v1.

\vspace{1ex}

\noindent [Tsu84]\, H. Tsuji --- {\it Complex Structures on $S^3\times S^3$} --- Tohoku Math. J (2) {\bf 36} (1984), no. 3, 351-376.

\vspace{6ex}

\noindent Dan Popovici, 

\noindent Institut de Math\'ematiques de Toulouse, Universit\'e Paul Sabatier, 

\noindent 118 Route de Narbonne, 31 062, Toulouse, France

\noindent Email: popovici@math.ups-tlse.fr

\end{document}